\documentclass[final,3p]{elsarticle}
\usepackage{}
\usepackage{amssymb}
\usepackage{graphicx}
\usepackage{amssymb}
\usepackage{amsmath}
\usepackage{cases}
\usepackage{amsfonts}
\usepackage{epsfig}
\usepackage{booktabs}
\usepackage{upgreek}
\usepackage{graphicx,tabularx,multirow, color}
\usepackage[active]{srcltx}
\usepackage{pifont}

\usepackage{float}%exact location
\usepackage{caption}
\usepackage{subcaption} % for subfigures
\usepackage[abs]{overpic}
\usepackage{url}

\usepackage[lined,commentsnumbered,ruled]{algorithm2e}

\makeatletter

\newcommand{\im}{\mathrm{i}}
\newcommand{\diff}{\mathrm{d}}
\newcommand{\rank}{\mathrm{rank}}
\newcommand{\diag}{\mathrm{diag}}
\newcommand{\spanset}{\mathrm{span}}

\newcommand{\tol}{\mathrm{tol}}
\newcommand{\gap}{\mathrm{gap}}

\newcommand{\cirrs}{RSRR}
\newcommand{\cirrm}{SSRR}

\newcommand{\bbx}{\mathbf{x}}
\newcommand{\bby}{\mathbf{y}}

\newtheorem{thm}{Theorem}

\makeatother

\begin{document}

\begin{frontmatter}



\title{Resolvent sampling based Rayleigh-Ritz method for large-scale nonlinear eigenvalue problems}
%A contour integral method for large-scale nonlinear eigenvalue problems with applications to BEM elastodynamics and acoustics


\author[a1]{Jinyou Xiao}
\ead{xiaojy@nwpu.edu.cn}
\cortext[cor1]{Corresponding author}
\address[a1]{School of Astronautics, Northwestern Polytechnical University, Xi'an 710072, China}

\author[a1]{Shuangshuang Meng}

\author[a3]{Chuanzeng Zhang\corref{cor1}}
\ead{c.zhang@uni-siegen.de}
\address[a3]{Department of Civil Engineering, University of Siegen, D-57068 Siegen, Germany}

\author[a4]{Changjun Zheng}
\ead{cjzheng@hfut.edu.cn}
\address[a4]{Institute of Sound and Vibration Research, Hefei University of Technology, Hefei, 230009, China}

\begin{abstract}

A new algorithm, denoted by \cirrs{}, is presented for solving large-scale nonlinear eigenvalue problems (NEPs) with a focus on improving the robustness and reliability of the solution, which is a challenging task in computational science and engineering. The proposed algorithm utilizes the Rayleigh-Ritz procedure to compute all eigenvalues and the corresponding eigenvectors lying within a given contour in the complex plane. The main novelties are the following. First and foremost, the approximate eigenspace is constructed by using the values of the resolvent at a series of sampling points on the contour, which effectively circumvented the unreliability of previous schemes that using high-order contour moments of the resolvent. Secondly, an improved Sakurai-Sugiura algorithm is proposed to solve the projected NEPs with enhancements on reliability and accuracy. The user-defined probing matrix in the original algorithm is avoided and the number of eigenvalues is determined automatically by provided strategies. Finally, by approximating the projected matrices with the Chebyshev interpolation technique, \cirrs{} is further extended to solve NEPs in the boundary element method, which is typically difficult due to the densely populated matrices and high computational costs. The good performance of \cirrs{} is demonstrated by a variety of benchmark examples and large-scale practical applications, with the degrees of freedom ranging from several hundred up to around one million. The algorithm is suitable for parallelization and easy to implement in conjunction with other programs and software.

\end{abstract}


\end{frontmatter}


\section{Introduction}
\label{s-intro}

The natural frequency $\lambda$ and natural mode $v$ are important parameters in the design of engineering structures and systems. Mathematically, they are the eigenvalue and the corresponding eigenvector satisfying the eigen-equation
\begin{equation}\label{gnep}
T(\lambda)v=0,
\end{equation}
where $T(z) \in \mathbb{C}^{n \times n}$ is a matrix-valued function depending on a parameter $z \in \mathcal{D} \subset \mathbb{C}$. When the matrix $T(z)$ in terms of $z$ is nonlinear, \eqref{gnep} becomes the so-called nonlinear eigenvalue problem (NEP). Except for the quadratic eigenvalue problem that has been well-studied in the dynamic analysis of structures \cite{TM01}, the numerical solution of general NEPs still remains a challenging task, especially at large scales \cite{MS11, effenberger2013robust, van2015rational}.

This work is concerned with developing efficient numerical methods for large-scale NEPs. In particular, we consider NEPs in the \emph{finite element method} (FEM) and the \emph{boundary element method} (BEM), since these two methods are extensively applied in computational science and engineering, either independently or in coupled manners \cite{OBZ15,Steinbach14K}. Their NEPs reflect the main bottlenecks in the current development of numerical methods for NEPs. In the FEM, NEPs are often caused by the $z$-dependent material properties and/or boundary conditions, examples including the analysis of structures involving viscoelastic materials \cite{DaPo01,CFEM06,ASPB09}, free vibration of plates and shells with elastically attached masses \cite{SSI06}, vibrations of fluid-solid structures \cite{CPV89}, special treatment of artificial boundary conditions \cite{BSS09}, etc. Typically, the FEM matrix $T(z)$ is sparse, symmetric, real and positive
or semi-positive definite. These properties greatly benefit the numerical treatment. However, the situations in the BEM are totally different. The eigenvalue problems are in general strongly nonlinear even if the underlying physical problem is linear \cite{KA93,ARY95,steinbach2012convergence}. The matrix $T(z)$ is typically complex, dense and unstructured, making the evaluation of $T(z)$ itself and the operations with $T(z)$ (e.g., applying to vectors, solving linear systems, etc) computationally very expensive. Besides, the entries of $T(z)$ are usually distinct functions of $z$ whose dependence on $z$ can hardly be expressed in explicit formulas. All these factors contribute to the immense difficulties of solving the NEPs in BEM even with the state-of-the-art techniques \cite{EK12, EC13}. Similar difficulties also exist in the NEPs in the coupled FEM-BEM \cite{EKSU12,Steinbach14K}.

In order to overcome the aforementioned difficulties in large-scale engineering applications, this paper aims at developing an eigensolver that simultaneously satisfies the following two requirements:
\begin{enumerate}[\it R1:]
  \item can robustly and reliably compute all the eigenvalues (and the corresponding eigenvectors) in a given region of the complex plane. \label{feature1}
  \item can solve large-scale NEPs (e.g., degrees of freedom (DOFs) $n$ up to several millions) no matter the matrix $T(z)$ is dense or sparse, structured or unstructured, and allows for an easy implementation and efficient parallelization. \label{feature2}
\end{enumerate}
We notice that although there exist a number of algorithms for the solution of NEPs (see e.g., \cite{MVVH04, Voss04, Voss07, Daniel09, Beyn12, EC13}), most of them are excluded by the above two requirements.

Essentially, existing numerical approaches for NEPs fall into three categories \cite{MVVH04, VRK13}. The first and classical approach is to formulate the NEP as a system of nonlinear equations and to use variations of Newton's method. Examples include the residual inverse iteration method \cite{Neu85}, the Jacobi-–Davidson method \cite{BV04, Voss07} and the block Newton method \cite{Daniel09}. Typically, Newton type methods are only locally convergent and need good initial guess of the eigenvalues. The Jacobi-–Davidson method can be used to determine several eigenvalues subsequently, but special strategies are required to avoid the repeated convergence towards the same eigenvalues \cite{Voss07, Daniel09}. Recently, a deflation strategy based on the concept of minimal invariant pairs has been proposed and integrated into the Jacobi-–Davidson method for the robust and successive computation of several eigenpairs \cite{EC13}. However, one should notice that, in all these methods, the repeated solution of the \emph{correction equation} during the iterative process can be computationally very expensive for large-scale problems. Newton type methods applied to small BEM eigenvalue problems can be found in, e.g., \cite{KAN96}.

The second approach transforms the NEP into a linear eigenvalue problem which can be solved by existing linear eigensolvers. This is a standard technique in solving quadratic eigenvalue problems in the FEM \cite{TM01}. The applications to the NEPs in the BEM can be found in \cite{ARY95}. Recently, linearization by the polynomial or rational approximation of $T(z)$ has attracted increasing attentions in solving medium and large NEPs; see, e.g., \cite{EK12, VRK13, KR14, guttel2014nleigs}. In particular, the compact rational Krylov method in \cite{van2015compact} exploits the structure of the linearization pencils by using a generalization of the compact Arnoldi decomposition. As a result, the extra memory and orthogonalization costs are negligible for large-scale problems.
Whereas, for the linearization-based methods there is still a need to construct linearizations that reflect the
structure of the given matrix polynomial and to improve the stability of the linearizations \cite{TM01,mackey2006structured}. Moreover, this class of methods are not suitable for large-scale BEM applications, because the storage of all the $n \times n$ coefficient matrices of $T(z)$ in a polynomial basis, even in compressed forms \cite{CWX15}, would require huge memory.

The third approach, called \emph{contour integral method}, is based on the contour integral of the resolvent $T(z)^{-1}$ applied to a full-rank matrix $U \in \mathbb{C}^{n \times L}$. This approach is first developed for solving generalized eigenvalue problems \cite{sakurai2003projection, polizzi2009density}, and later, extended to solving NEPs in \cite{SAT09} and \cite{Beyn12}, respectively, using the Smith form and Keldysh's theorem for analytic matrix-valued functions. Since \cite{SAT09} and \cite{Beyn12} come up with similar algorithms, we consider the block Sakurai-Sugiura (SS) algorithm proposed in \cite{SAT09}. This algorithm transforms the original NEP into a small-sized generalized eigenvalue problem involving block Hankel matrices by using the monomial moments of the probed resolvent $T(z)^{-1} U$ on a given contour. Applications of the block SS algorithm for solving NEPs in the BEM are presented in \cite{steinbach2012convergence, Gao13, LL13, ZCJ15}. The algorithm could be efficient in solving general NEPs only under the following conditions: (1) the number of columns of $U$, i.e., $L$, is sufficiently large, and (2) the contour integrals are accurately evaluated. Otherwise, it would become unreliable and inaccurate \cite{sakurai2007cirr, YS13}, and loss of eigenvalues and spurious eigenvalues may frequently occur. But meeting these conditions, if possible, often implies heavy computational burdens in large-scale applications.

Although still under
development, the contour integral method appears to be a very promising candidate that potentially meets the above two requirements \emph{R\ref{feature1}} and \emph{R\ref{feature2}}. It simultaneously considers all the eigenvalues within a given contour; the most computationally intensive part, i.e., the solution of a series of independent linear systems, can be effectively parallelized. The main difficulty is to enhance the robustness and accuracy while still keeping the computational cost to a reasonable level. Along this line, the \underline{S}akurai-\underline{S}ugiura method with the \underline{R}ayleigh-\underline{R}itz projection (\cirrm{}) has been recently proposed \cite{YS13}. In \cirrm{} the eigenspaces for projection are generated from a matrix $\hat M \in \mathbb{C}^{n \times K\cdot L}$ collecting all the monomial moments of the resolvent $T(z)^{-1} U$ up to order $K$. We refer to this as the \emph{resolvent moment scheme} (or simply, \emph{moment scheme}) for generating eigenspaces. Numerical results have shown that
the \cirrm{} algorithm can significantly improve the robustness and accuracy compared with the block SS algorithm.

Despite the great potential, the \cirrm{} algorithm suffers from the possible failure in practical applications. The reason is indeed easy to understand: the eigenspaces generated by the moment scheme, $\hat M$, tends to be unreliable when high order $K$ of moments are involved, as the matrix $\hat M$ would become rank-deficient and consequently the rank condition $\rank (\hat M) \geqslant \bar n_{\mathcal {C}}$ for a proper extraction of all the $\bar n_{\mathcal {C}}$ eigenpairs inside the contour $\mathcal {C}$ will never be satisfied; see Section \ref{S-eigenspace-mom} for detailed reasoning, and Sections \ref{S-ne-4bench} to \ref{S-ne-bem} for a variety of numerical evidences. Note that, besides increasing $K$, the condition $\rank (\hat M) \geqslant \bar n_{\mathcal {C}}$ can also be fulfilled by increasing the value of $L$. But bear in mind that the main computational work of the \cirrm{}, i.e., the computation of $T(z)^{-1} U$, is in proportion to $L$, thus a small $L$ is essential for reducing the total computational costs in large-scale solutions.

To settle the possible failure in eigenspace construction, we propose a more reliable \emph{resolvent sampling scheme} (or simply, \emph{sampling scheme}). In this scheme, the eigenspaces are generated by a matrix $\hat S \in \mathbb{C}^{n \times NL}$ collecting all the matrices $T(z_k)^{-1} U$ at a series of sampling points $z_k \,( k=1, \cdots, N)$, instead of the moments of $T(z)^{-1} U$.
The sampling points can either be a set of quadrature nodes on the boundary of the domain of interest, or a general set of points in the domain.
The sampling scheme has several advantages over the moment scheme: (1) the possible failure in eigenspace construction caused by the use of high order moments is effectively remedied and the algorithm is more concise; and more improtantly, (2) the eigenspaces are better than those generated by the moment scheme and thus the eigensolutions are more accurate. Besides, the sampling scheme inherits the salient features of the contour integral method, i.e., simultaneous consideration of all eigenvalues within a given contour and good parallelizability.

We also propose a \emph{\underline{r}esolvent \underline{s}ampling based \underline{R}ayleigh-\underline{R}itz method}, abbreviated as \cirrs{}, for solving general NEPs. The projected NEPs of small sizes are solved using the block SS algorithm \cite{SAT09}. As mentioned before, the block SS algorithm can be robust and accurate under two conditions, which can be easily satisfied for the reduced NEPs since they are often of very small sizes (e.g., several tens or hundreds).
Also due to the small size of the object problems, we directly compute the moments of the resolvent of the reduced matrix instead of using the probing matrix $U$ in the block SS algorithm. As such, only very low order moments are needed, greatly enhancing the robustness of the algorithm. Moreover, a strategy for automatically determining the number of the eigenvalues is also proposed.

The \cirrs{} algorithm is able to meet the previously mentioned two requirements \emph{R\ref{feature1}} and \emph{R\ref{feature2}}. It can be used to solve general NEPs no matter the matrix $T(z)$ is sparse or dense, structured or unstructured, provided that a fast linear system solver for computing $T(z_k)^{-1} U$ and an efficient routine for computing the reduced matrix of $T(z)$ are available. It permits an efficient parallelization and can be easily implemented in conjunction with other programs or software. For example, by combining the \cirrs{} with the FEM software \textsc{Ansys}\textregistered, we can solve FEM NEPs with more than 1 million DOFs on a Server; see examples in Section \ref{S-ne-fem}. The implementation of the \cirrs{} in the BEM is not so straightforward as in the FEM, since an explicit expression of $T(z)$ in terms of $z$ is in general not available for the BEM. We resolve this problem by computing the Chebyshev polynomial interpolation of $T(z)$. The price to pay is a moderate increase of time, but the operation is very suitable for parallelization. Finally, we are able to solve BEM NEPs with several tens of thousands of DOFs on a personal computer; see Section \ref{S-ne-bem}.

Finally, we notice that there are several methods that are specially designed for solving NEPs in the FEM. Examples include the order-reduction-iteration method and the asymptotic method; see \cite{DaPo01, CFEM06, bilasse2009generic} and the references therein. These methods typically solve a simplified problem, for example, the eigenvalues and associated eigenvectors in the desired region for undamped problems are used for the projection. They are partially heuristic since there is no guarantee that all the interested eigenvalues are captured or that the eigensolutions of the projected system are close to those of the original system \cite{MS11}.

The rest of this paper is organized as follows.
In Section \ref{S-basics}, the basic Keldysh's theorem and the Rayleigh-Ritz procedure are briefly reviewed.
In Section \ref{S-eigenspace}, the moment scheme and sampling scheme for constructing eigenspaces are described.
 By using the Rayleigh-Ritz projection, the NEP \eqref{gnep} is transformed into a reduced NEP of a small size. A modified block SS algorithm is proposed for the robust and accurate solution of the reduced NEP in Section \ref{S-reducedNEP}.
 The \cirrs{} algorithm is summarized in Section \ref{S-algorithm}, and is extended to solve the NEPs in the acoustic BEM in Section \ref{S-bem}.
In Section \ref{S-ne} the performance of the \cirrs{} is attested by a variety of representative examples and comparisons. Section \ref{S-conclusion} presents the essential conclusions of the paper.

\section{Mathematical foundations} \label{S-basics}

Consider a holomorphic matrix-valued function $T(z) \in \mathbb{C}^{n \times n}$ defined in an open domain $\mathcal {D} \in \mathbb{C}$, and assume that $T(z)$ is regular in the sense that the determinant, $\det T(z)$, does not vanish identically. We intend to search for all the eigenvalues within a compact set $C \subset \mathcal {D}$ and the associated eigenvectors. Our method is motivated by the following theorem regarding the relation between the eigenvalues and the resolvent $T(z)^{-1}$ (see \cite{Beyn12}, Corollary 2.8).
\begin{thm} \label{thm1}
Let $C \subset \mathcal {D}$ be a compact set containing a finite number $n_C$ of different eigenvalues $\lambda_k\, (k=1, \cdots, n_C)$, and let
\begin{equation} \label{eq-VC}
V_C = \left ( v^{l,k}_j, \quad 0 \le j \le \mu_{l,k} -1, \, 1 \le l \le \eta_k, \, k=1, \cdots, n_C \right)
\end{equation}
and
\begin{equation} \label{eq-WC}
W_C = \left ( w^{l,k}_j, \quad 0 \le j \le \mu_{l,k} -1, \, 1 \le l \le \eta_k, \, k=1, \cdots, n_C  \right)
\end{equation}
be the corresponding \emph{canonical systems of generalized eigenvectors} (CSGEs) of $T$ and $T^H$, respectively. Then there exists a neighborhood $C \subset \mathcal {U} \subset \mathcal {D}$ and a holomorphic matrix-valued function $R: \Omega \rightarrow \mathbb{C}^{n \times n}$ such that for all $z \in \mathcal {U} \setminus \{\lambda_1, \cdots, \lambda_{n_C }\}$,
\begin{equation}\label{eq-thm1}
T(z)^{-1} = \sum^{n_{C}}_{k=1}  \sum^{\eta_k}_{l=1}  \sum^{\mu_{l,k}}_{j=1}  ( z-\lambda_k )^{-j}  \sum^{\mu_{l,k}-j}_{m=0}    v^{l,k}_m \left(w^{l,k}_{\mu_{l,k} -j-m}\right)^{H} + R(z).
\end{equation}
\end{thm}

In expressions \eqref{eq-VC}, \eqref{eq-WC} and \eqref{eq-thm1}, $\eta_k$ and $\mu_{l,k}$ represent the dimension of the nullspace of $T(\lambda_k)$ and the $l$th partial multiplicity of $T(z)$ at $\lambda_k$, respectively. We refer the readers to \cite{Beyn12, YS13} for more detailed description of the related notations. The superscript $H$ denotes the conjugate transpose.

Theorem \ref{thm1} shows that the eigenvalues of $T(z)$ are the poles of $T(z)^{-1}$. By means of the residue theorem, it is possible to extract the eigenpairs inside a given contour, denoted by $\mathcal{C}$, by performing contour integration on $\mathcal {C}$. This fact leads to the block SS algorithm \cite{AS09} and Beyn's algorithm \cite{Beyn12} for solving NEPs. However, as previously stated, this class of algorithms may suffer from instability issues and high computational burden in large-scale applications. These problems can be effectively circumvented by using the Rayleigh-Ritz procedure.

The classical Rayleigh-Ritz procedure relies on a proper search space $\mathcal{S}$ that contains the interested eigenvalues. Once an orthogonal basis $S\in \mathbb{C}^{n \times k_S}$ of the search space $\mathcal{S}$ is obtained, the original NEP \eqref{gnep} can be transformed to the reduced NEP
\begin{equation} \label{eq-Ts-NEP}
T_S(\lambda)g = 0 \quad \mathrm{with} \quad T_S(z) = S^H T(z)S \in \mathbb{C}^{k_S \times k_S}.
\end{equation}
Let $(\lambda, g)$ be any eigenpair of the reduced NEP, then $(\lambda, Sg)$ is an eigenpair of the original NEP \cite{YS13}.

Therefore, the construction of appropriate search spaces is of primary importance to the Rayleigh-Ritz procedure. In the following Section \ref{S-eigenspace} we will propose a new scheme for this purpose, and then in Section \ref{S-reducedNEP} we will develop a modified block SS algorithm for solving the reduced NEP \eqref{eq-Ts-NEP}.

\section{Generation of approximate eigenspaces} \label{S-eigenspace}

\begin{figure}[hbt]
\centering
\epsfig{figure=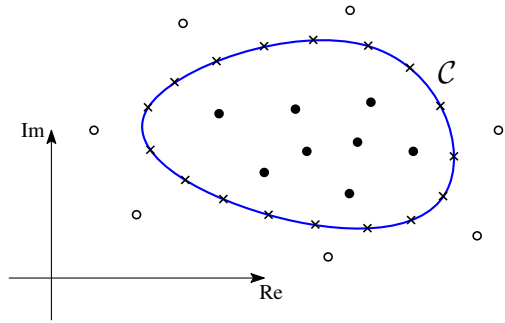,width=0.4\textwidth}
\caption{A diagram showing the contour $\mathcal {C}$ (solid line), eigenvalues interested ($\bullet$), eigenvalues not interested ($\circ$) and quadrature (or sampling) points on the contour ($\scriptstyle\pmb{\times}$).}
\label{fig-contour}
\end{figure}

Now consider a contour $\mathcal {C} \subset \mathcal {D}$ and assume that $T(z)$ has no eigenvalues on the contour $\mathcal {C}$; see Figure \ref{fig-contour}. Denote by $n_{\mathcal {C}}$ the number of different eigenvalues within $\mathcal {C}$, and by $\bar n_{\mathcal {C}}$ the total number of eigenvalues, counting the algebraic multiplicity, i.e.,
$\bar n_{\mathcal{C}} = \sum^{n_{\mathcal{C}}}_{k=1}  \sum^{\eta_k}_{l=1} \mu_{l,k}$. Let $L$ be an integer larger than or equal to the maximal algebraic multiplicity of the eigenvalues within $\mathcal {C}$, which will be used to denote the number of vectors used to probe the resolvent $T(z)^{-1}$. Usually $L$ is a small number independent of $n$.

We intend to compute all the $\bar n_{\mathcal {C}}$ eigenpairs within $\mathcal {C}$ using the Rayleigh-Ritz procedure.
In the following, we first briefly describe the moment scheme for generating the projection subspace $\mathcal{S} \approx \spanset(V_\mathcal {C})$ and show the possible failure of the scheme.
Then, we propose the new sampling scheme and show that it can always generate more reliable eigenspaces.

\subsection{Resolvent moment scheme and its failure} \label{S-eigenspace-mom}

Given a matrix $U \in \mathbb{C}^{n \times L}$ with $L$ linear independent columns, the moment matrices of the probed resolvent $T(z)^{-1}U$ can be defined as
\begin{equation}\label{eq-malpha}
M_{\alpha} = {1 \over 2 \pi \im} \int_\mathcal{C} z^\alpha  T(z)^{-1} U \diff z = V_\mathcal{C}  \Lambda^\alpha  W_\mathcal{C}^H, \quad \alpha = 0, \cdots, K'-1,
\end{equation}
where $V_\mathcal{C}, W_\mathcal{C} \in \mathbb{C}^{n \times \bar n_{\mathcal{C}}}$ are the matrices consist of the CSGEs of $T$ and $T^H$, respectively; $\Lambda$ has Jordan normal form given by
\begin{equation} \label{eq-Mpr-Jordan}
\Lambda = \begin{bmatrix}
  J_1 &   &   \\
    & \ddots &   \\
    &   & J_{n_\mathcal {C}} \\
\end{bmatrix}, \quad
J_k = \begin{bmatrix}
  J_{k,1} &   &   \\
    & \ddots &   \\
    &   & J_{k,\eta_k} \\
\end{bmatrix}, \quad
J_{k, l} \equiv J_{\lambda_k, \mu_{l,k} } = \begin{bmatrix}
           \lambda_k  & 1 &   &   \\
             & \ddots & \ddots &   \\
             &   & \lambda_k & 1 \\
             &   &   & \lambda_k  \\
         \end{bmatrix} \in \mathbb{C}^{\mu_{l,k} \times \mu_{l,k}}.
\end{equation}
The second equality of \eqref{eq-malpha} is obtained by inserting \eqref{thm1} into the contour integral and using the Cauchy residual theorem; see equation (49) in \cite{Beyn12}.
Note that, besides the monomial functions $z^\alpha$, other basis functions can be used as well.

By using the relation $M_{\alpha} = V_\mathcal{C}  \Lambda^\alpha  W_\mathcal{C}^H$, one has
\begin{equation}\label{eq-Mrow}
M := \left ( M_0, \, M_1, \, \cdots, \, M_{K'-1} \right) = V_{\mathcal{C}} W_{\mathcal{C}}^{[K']}
\end{equation}
where
$$
W_{\mathcal{C}}^{[K']} = \begin{bmatrix}
      W^H_{\mathcal{C}} U & \Lambda W^H_{\mathcal{C}} U & \cdots & \Lambda^{K'-1} W^H_{\mathcal{C}} U  \\
    \end{bmatrix} \in \mathbb{C}^{ \bar n_{\mathcal{C}} \times K'\cdot L}.
$$
In order to retrieve the eigenspace $\spanset (V_{\mathcal{C}})$ from the matrix $M$, the numbers $K'$ and $L$ have to be chosen such that
\begin{equation} \label{eq-rank-M-KL}
K'L \ge \rank( M) =  \rank(V_{\mathcal{C}}).
\end{equation}

From the viewpoint of reducing the computational costs, a small $L$ but a large $K'$ are preferable \cite{YS13} to ensure the rank condition \eqref{eq-rank-M-KL}. However, we notice that the quality of the eigenspace constructed by the moment scheme
would deteriorate when a large $K'$ is used, because, in finite precision arithmetic, the columns of $M$ tend to be linearly dependent when $K'$ is large, and finally the condition \eqref{eq-rank-M-KL} may never be fulfilled accurately. As a result, some computed eigenpairs might be of low accuracy, or more annoyingly, spurious eigenvalues and/or loss of eigenvalues might occur.

As an illustration, we assume that all eigenvalues are simple so that $L=1$ can be used. In this case $n_{\mathcal{C}} = \bar n_{\mathcal{C}}$. Then, similar to \eqref{eq-Mrow},
$M$ can be written as
$$
M = V_{\mathcal{C}} D B,
$$
where $D\in \mathbb{C}^{n_{\mathcal{C}} \times n_{\mathcal{C}}}$ is a diagonal matrix whose diagonal entries are given by the column vector $W_{\mathcal{C}}^H U $, and
$$
B = \begin{bmatrix}
      1                & \lambda_1                   & \cdots & \lambda_1^{K'-1} \\
      1                & \lambda_2                   & \cdots & \lambda_2^{K'-1} \\
      \vdots           & \vdots                      & \ddots &          \vdots \\
      1                & \lambda_{n_{\mathcal{C}}}   & \cdots & \lambda_{n_{\mathcal{C}}}^{K'-1} \\
    \end{bmatrix}.
$$

\begin{figure}[hbt]
\centering
\epsfig{figure=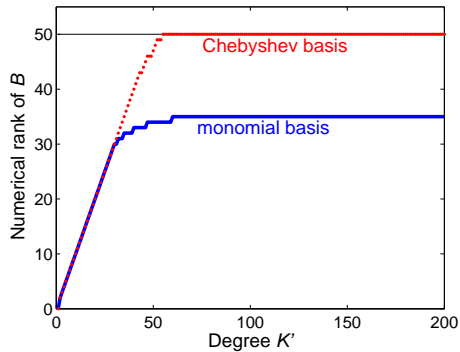,width=0.4\textwidth}
\caption{Numerical rank of $B$ versus $K'$ in the monomial basis and Chebyshev basis.}
\label{fig-rankB}
\end{figure}

In order to extract all the $n_{\mathcal{C}}$ eigenvalues, the numerical rank of $B$, denoted by $k_B$, has to be equal to $n_{\mathcal{C}}$. However, typically this condition is hard to be satisfied for relatively large $n_{\mathcal{C}}$, due to the fact that $B$ is a Vandermonde matrix which becomes extremely ill-conditioned in large dimensions.
For example, let $n_{\mathcal{C}} = 50$, then condition \eqref{eq-rank-M-KL} implies that $k_B$ has to be at least 50. Let us further assume that the eigenvalues are equally distributed in the interval $[-0.9,0.9]$, which can be computed by using Matlab function \verb"linspace"$(-0.9,0.9,50)$. Let $\sigma_1 \ge \sigma_2 \ge \cdots \ge \sigma_{50}$ be the singular values of $B$. The ranks $k_B$ is determined by condition that $\sigma_{k_B}\ge tol\cdot\sigma_1$ but $\sigma_{k_B+1} < tol\cdot\sigma_1$, where we set $tol=10^{-12}$. The relation between $k_B$ and $K'$ is shown by the blue line in Figure \ref{fig-rankB}. Obviously, $k_B$ never exceeds 35 no matter how large $K'$ is!

The example just given is admittedly manufactured; however, it is easy to imagine that similar situations could arise in practical applications (see Section \ref{S-ne} for several examples).
One remedy for the failure in the eigenspace generation is to choose a right set of basis functions for the moments \eqref{eq-malpha}. For example, if we use Chebyshev functions of the first kind instead of the monomial functions, the condition $k_B = n_{\mathcal{C}}=50$ can be fulfilled when $K' \ge 55$, as indicated by the red line in Figure \ref{fig-rankB}. However, in the following section we will propose a more general and reliable remedy in which the choices of basis functions and parameter $K'$ can be completely avoided.

\subsection{Resolvent sampling scheme} \label{S-eigenspace-samp}

We begin by noting that, when using the moment scheme, the contour integral in $M_\alpha$ \eqref{eq-malpha} has to be evaluated by numerical quadratures, leading to the following approximate moments
\begin{equation}\label{eq-malpha-qud}
\hat M_\alpha := \sum_{i=0}^{N-1} \omega_i z_i^\alpha T(z_i)^{-1} U \approx M_{\alpha}, \quad \alpha=0,1,\cdots, K'-1,
\end{equation}
where $N$ is the number of quadrature nodes, $z_i$ and $\omega_i \, (i=0, \cdots, N-1)$ are quadrature nodes and weights, respectively. Here, $N$ is assumed to be large enough so that the quadrature error is negligible. Then the moment scheme actually generates eigenspaces by
\begin{equation}\label{eq-hatM}
\hat M = \left ( \hat M_0, \, \hat M_1, \, \cdots, \, \hat M_{K'-1} \right) \in \mathbb{C}^{n \times K'\cdot L }.
\end{equation}
Similar to its analytical counterpart $M$, the matrix $\hat M$ would become rank-deficient for large $K'$.

To motivate the sampling scheme, we further exploit the structure of the matrix $\hat M$. From \eqref{eq-malpha-qud} it can be seen that each $\hat M_\alpha$
is indeed a linear combination of $T(z_i)^{-1} U$. Thus, by letting $\hat S$ to be the \emph{sampling matrix} collecting all the probed resolvent $T(z_i)^{-1} U$,
\begin{equation} \label{eq-hatS}
\hat S = \left [ T(z_0)^{-1} U, \, T(z_1)^{-1} U, \, \cdots, \, T(z_{N-1})^{-1} U \right] \in \mathbb{C}^{n \times N\cdot L},
\end{equation}
the matrix $\hat M$ can then be expressed as
\begin{equation} \label{eq-hatMS}
\hat M = \hat S D Z,
\end{equation}
where $D \in \mathbb{C}^{N\cdot L \times N\cdot L}$ is a diagonal matrix, $D = \diag ( D_i,\, i = 0,1,\cdots, N-1 )$, with $D_i$ being $L\times L$ diagonal matrices whose diagonal entries are all $\omega_i$; and $Z$ is a Vandermonde-like matrix
$$
Z = \begin{bmatrix}
           Z_0^{[0]}       & Z_0^{[1]}         & \cdots  &  Z_0^{[K'-1]} \\
           Z_1^{[0]}       & Z_1^{[1]}         & \cdots  &  Z_1^{[K'-1]} \\
                     \vdots      & \vdots  & \ddots  & \vdots \\
           Z_{N-1}^{[0]}       & Z_{N-1}^{[1]}     & \cdots  &  Z_{N-1}^{[K'-1]}  \\
         \end{bmatrix} \in \mathbb{C}^{N\cdot L \times K'\cdot L},
$$
with each block $Z_i^{[k]}$ being a $L\times L$ matrix whose entries are all $z_i^{k}$. Obviously, $Z$ would become rank-deficient when $K'$ is large, which accounts for the possible rank-deficiency of $\hat M$. More importantly, \eqref{eq-hatMS} implies that the rank property of $\hat S$ should be better than that of $\hat M$.

Another important implication of \eqref{eq-hatMS} is that the moment matrix $\hat M$ lies in the range of the sampling matrix $\hat S$,
\begin{equation} \label{eq-setsMS}
\spanset (\hat M) \subseteq \spanset (\hat S),
\end{equation}
which inspires us to construct eigenspaces directly using the sampling matrix $\hat S$, instead of using the
moment matrix $\hat M$. That is why it is called ``the sampling scheme''.

The sampling scheme indeed constructs a better subspace for projection. On one hand, the two schemes can obtain the same subspaces only if $\hat M$ is not rank-deficient and $K'=N$. Otherwise, in most situations the subspace of the sampling scheme is always larger than that of the moment scheme. Further, relation \eqref{eq-hatMS} indicates that the numerical rank of $\hat S$ should not be smaller than that of $\hat M$, which means that the subspace of the sampling scheme could be more stable for projection. This argument will be largely validated by the numerical examples in Section \ref{S-ne}.

On the other hand, the sampling scheme can generate a reliable subspace that contains all the eigenvalues, provided that $N$ is large enough; that is,
\begin{equation} \label{eq-VCS}
      \spanset (V_{\mathcal{C}}) \subseteq \spanset (\hat S), \quad N\rightarrow \infty.
      \end{equation}
 Relation \eqref{eq-VCS} can be verified by using the relation \eqref{eq-setsMS} and Theorem 3 in \cite{YS13}. In fact, the latter shows that $\hat M$ tends to be equal to $M$ \eqref{eq-Mrow} when $N$ increases and the quadrature rule in \eqref{eq-malpha-qud} converges. Hence, there exists a positive integer $K_0$ such that
      \begin{equation} \label{eq-VChatM}
      K_0< K'\le N \quad \Rightarrow \quad \spanset (V_{\mathcal{C}}) \subseteq \spanset (\hat M).
      \end{equation}
      Note that the following condition is a prerequisite for \eqref{eq-VChatM}
      \begin{equation} \label{eq-rank-hatM}
      K'\cdot L \ge  \rank(\hat M) \ge  \rank(V_{\mathcal{C}}).
      \end{equation}
      The above arguments also indicate that for constructing effective subspaces that contain all the eigenvalues, $N$ has to satisfy the following condition
      \begin{equation} \label{eq-rank-NL}
      N\cdot L \ge  \rank(\hat S) \ge \rank(V_{\mathcal{C}}).
     \end{equation}

In the Rayleigh-Ritz procedure, an orthogonal basis of $\hat S$, denoted by $S$, is needed. It can be computed by the truncated singular value decomposition (SVD) of $\hat S$ with a tolerance $\delta$.
Since the number of the columns of $\hat S$ and $N\cdot L$, is typically small and independent of the size $n$ of the problem, the computational cost of the SVD scales only linearly with $n$.
Let $k_S$ denote the numerical rank of $\hat S$.
In order to properly extract all the eigenvalues inside $\mathcal {C}$, the condition $k_S \ge \rank(V_{\mathcal{C}})$ has to be satisfied. However, since $\bar n_{\mathcal {C}} \ge \rank(V_{\mathcal{C}})$, it is more convenient to use the following condition in practice
\begin{equation} \label{eq-rank-nc}
N\cdot L \ge k_S \ge \bar n_{\mathcal {C}},
\end{equation}
which means that $N$ and $L$ have to be chosen such that $k_S$ is larger than the total number of the eigenvalues.
Of course, $L$ should not be less that the maximal algebraic multiplicity of the eigenvalues in $\mathcal {C}$.

We close this section by several comments on the connections and comparisons between the proposed sampling scheme and some existing nonlinear eigensolvers.

First, although the sampling scheme is motivated by the contour integral based moment scheme, it is more general than the latter. The connection of the two schemes lies in the fact that both of them construct eigenspaces based on the solutions $Y_i$ of a series of independent linear system of equations $T(z_i)Y_i= U, \, i=0, \cdots, N-1$, and the subspace of the moment scheme converges to that of the sampling scheme when $K'\rightarrow N$ and the rank deficiency does not occur. However, the contour moment scheme is restricted to simple closed contours in the complex plane, while there is nothing in the sampling matrix \eqref{eq-hatS} that constrains the sampling points $z_i$ to be chosen according to a quadrature rule on a certain contour. Indeed, our numerical experience shows that the sampling points within the contour can also be used in the sampling scheme, which sometimes may even improve the accuracy of the eigen-solutions. Whereas, in this paper we shall not further pursue this topic, but always choose $z_i$ on the contour since this choice can lead to satisfactory results in general situations.

%Second, the moment scheme has an advantage over the sampling scheme: the number of columns of $\hat M$ is often less than that of $\hat S$. As a result, the cost for computing orthogonal bases for the eigenspaces is much less in the moment scheme, since this part of computational cost scales quadratically with the number of the columns.
%This advantage could be remarkable when a large number (e.g., hundreds) of eigenvalues are sought simultaneously. In this situation, however, the computational cost of the sampling scheme can be reduced by a subdivision of the domain of eigenvalues and using it for each sub-domain.

Second, the sampling scheme has connections with several other nonlinear eigensolvers. As previously stated, the essence of the sampling scheme is to construct adequate eigenspaces using the solutions of the equations $T(z_i)Y_i= U$ at a series of sampling points $z_i$.
Indeed, there exist several other eigensolvers that are based on solving the equations with the matrices $T(z_i)$, but their sampling points $z_i$ and right-hand-sides (rhs) are generally obtained in different manners. For example, in the classical inverse iteration and Rayleigh functional iteration, the rhs are determined by the previous iterations and the derivatives of the matrix $T(z)$ \cite{MVVH04}. Recently, a new class of nonlinear eigensolvers based on the dynamic linearization and the rational Krylov method have been proposed; see \cite{guttel2014nleigs} and the references therein. A representative example is the NLEIGS method in \cite{guttel2014nleigs}. The NLEIGS utilizes a dynamically structured rational interpolant of $T(z)$ and a new companion-type linearization to transform the NEP into a generalized eigenvalue problem with a special structure, which is then solved using the rational Krylov method. Due to the special structure, only the equations with the matrices $T(z_i)$ have to be solved at each iteration of the rational Krylov method, instead of the much larger equations with the linear pencils; see Lemma 4.5 and Corollary 4.6 in \cite{guttel2014nleigs}. This class of methods may have the advantage of using a smaller number of sampling points $z_i$ due to the endeavors in their optimization, but the rhs of the equations are also determined by the previous iterations. On the contrary, all the equations in the sampling scheme have the same rhs $U$, which is very suitable for parallelization.
A numerical comparison of the present \cirrs{} and the NLEIGS is conducted in Section \ref{S-ne-comparison}.

\section{Solution of the reduced NEP} \label{S-reducedNEP}

Once the eigenspaces are obtained, by using the standard Rayleigh-Ritz procedure the original NEP \eqref{gnep} is transformed into the reduced NEP \eqref{eq-Ts-NEP}.
The dimension of $T_S(z)$, $k_S$, usually ranges from several tens to several hundreds. Quite a few existing methods can be used to handle NEPs of this type. Whereas, a method that can robustly and accurately capture all eigenvalues in certain region is still lacking. The block SS algorithm proposed in \cite{AS09} (see also \cite{Beyn12}) is capable of directly computing all the eigenvalues within a given contour, without needing an initial guess of their locations. But there are two important parameters left for the users to determine: one, denoted by $L'$ \footnote{Notice that the parameter $L'$ here is denoted by $L$ in \cite{AS09}. It plays a similar role as the parameter $L$ in Section \ref{S-eigenspace} for generating eigenspaces. }, has to be larger than the maximal algebraic multiplicity of the eigenvalues inside $\mathcal {C}$; another is the threshold used to truncate the singular values of the Hankel matrix $H$. Both the two parameters are problem dependent. Improper values of them can lead to wrong eigenvalues.

%For instance, in locally convergent methods of Arnoldi or Jacobi-Davidson type (see \cite{Voss04,BV04}), an initial guess of eigenvalues or eigenspaces is typically needed. Even though the convergence of these method can be problematic when the so-called minimum-maximum characterizations are not available \cite{Daniel09,Beyn12}.

Here, we modify the block SS algorithm so that the parameter $L'$ is avoided. Meanwhile, we put forward a strategy to automatically determine the number of singular values of $H$ that should be retained in the truncation. The basic theory involved here can be found in \cite{AS09,Beyn12}.

\subsection{Moments of $T_S(z)^{-1}$   }

The block SS algorithm starts by computing the following moments of the matrix function $T_S(z)^{-1}$ up to order $2K-1$
\begin{equation} \label{eq-Ap}
\begin{aligned}
A_\alpha &= {1 \over 2 \pi \im} \int_{\mathcal {C}} \left({z-\gamma \over \rho }\right)^\alpha \, T_S(z)^{-1}  \diff z, \quad \alpha=0,1,\cdots, 2K-1,
\end{aligned}
\end{equation}
where $K$ denotes the number of blocks in each row or column of the Hankel matrices; see \eqref{eq-HH}. Note that: (i) here we use the same contour $\mathcal {C}$ as the one used in generating the eigenspace in section \ref{S-eigenspace}; (ii) the shift and scaling parameters $\gamma$ and $\rho$ are introduced to improve the stability, and for the elliptical contour defined in \eqref{eq-contour} $\gamma$ is the center of the ellipse and $\rho = \max(a,b)$; (iii) the moments are different from those used in \cite{AS09} since the random probing matrix of dimension $n \times L'$ is not explicitly involved here due to the small size of the problem. From Theorem \ref{thm1} one has,
\begin{equation} \label{eq-Adecomp}
A_\alpha = \tilde V \Lambda^\alpha \tilde W^H,
\end{equation}
where
$$
\begin{aligned}
\tilde V = \left ( \tilde v^{l,k}_j, \quad 0 \le j \le \mu_{l,k} -1, \, 1 \le l \le \eta_k, \, 1 \le k \le n_\mathcal {C} \right), \\
\tilde W = \left ( \tilde w^{l,k}_j, \quad 0 \le j \le \mu_{l,k} -1, \, 1 \le l \le \eta_k, \, 1 \le k \le n_\mathcal {C} \right)
\end{aligned}
$$
are the CSGEs of $T_S$ and $T_S^H$ defined analogously as \eqref{eq-VC} and \eqref{eq-WC}. The matrix $\Lambda$ has Jordan normal form similar to that in \eqref{eq-Mpr-Jordan}, but with Jordan blocks $J_{k,l} = J_{\lambda'_k, \mu_{l,k} }$,
$\lambda'_k = {\lambda_k - \gamma \over \rho}$.

\subsection{ Transformed linear eigenvalue problem and its solution}
 We define two block Hankel matrices as
\begin{equation} \label{eq-HH}
H= \begin{bmatrix}
      A_0     & A_1 & \cdots & A_{K-1} \\
      A_1     & \ddots     & \ddots &         \vdots \\
      \vdots         & \ddots     & \ddots &          \vdots \\
      A_{K-1} & \cdots     & \cdots & A_{2K-2} \\
    \end{bmatrix}, \quad
H^{<}= \begin{bmatrix}
      A_1       & A_2 & \cdots & A_{K} \\
      A_2       & \ddots     & \ddots &           \vdots \\
      \vdots           & \ddots     & \ddots &           \vdots \\
      A_{K} & \cdots     & \cdots & A_{2K-1} \\
    \end{bmatrix}.
 \end{equation}
By using \eqref{eq-Adecomp}, one has
\begin{equation*} \label{H-ss}
H= \begin{bmatrix}
      \tilde V \\
      \tilde V\Lambda \\
      \vdots \\
      \tilde V\Lambda^{K-1} \\
    \end{bmatrix}
    \begin{bmatrix}
      \tilde  W^H & \Lambda \tilde  W^H & \cdots & \Lambda^{K-1} \tilde  W^H \\
    \end{bmatrix}, \quad
H^{<}= \begin{bmatrix}
      \tilde V \\
      \tilde V\Lambda \\
      \vdots \\
      \tilde V\Lambda^{K-1} \\
    \end{bmatrix}
    \Lambda
    \begin{bmatrix}
      \tilde  W^H & \Lambda \tilde  W^H & \cdots & \Lambda^{K-1} \tilde  W^H \\
    \end{bmatrix}.
 \end{equation*}
To ensure that the eigenvalue matrix $\Lambda$ can be properly extracted from the matrices $H$ and $H^<$,
$K$ has to be chosen such that it holds
\begin{equation} \label{eq-ranks}
\rank \begin{bmatrix}
      \tilde V \\
      \tilde V\Lambda \\
      \vdots \\
      \tilde V\Lambda^{K-1} \\
    \end{bmatrix}
    =
\rank \begin{bmatrix}
      \tilde  W^H & \Lambda \tilde  W^H & \cdots & \Lambda^{K-1} \tilde  W^H \\
    \end{bmatrix}
    = \bar n_{\mathcal {C}}.
\end{equation}
In most real applications, a very small $K$ is often sufficient, e.g., $K < 10$.

 To extract the eigenvalues, one first compute the truncated SVD of $H$,
\begin{equation} \label{eq-svdH}
 H = V \Sigma W^H \approx V_0 \Sigma_0 W_0^H.
\end{equation}
The way to obtain $V_0$, $\Sigma_0$  and $W_0$ from $V$ , $\Sigma$ and $W$ will be detailed in Section \ref{S-S-ncs}. Then we compute the matrix
\begin{equation} \label{eq-A}
A = V_0^H H^{<} W_0 \Sigma_0^{-1},
\end{equation}
which has the Jordan normal form $\Lambda$ \cite{Beyn12}. The eigenvalue problem for $A$ can be easily solved using the Matlab function \verb"eig". Let $(\lambda', g')$ be an eigenpair of $A$, then the corresponding eigenpair $(\lambda, g)$ of the reduced NEP \eqref{eq-Ts-NEP} can be obtained as
\begin{equation} \label{eq-eigens}
g = H_r W_0 \Sigma_0^{-1} g', \quad \lambda = \rho\lambda' + \gamma,
\end{equation}
where $H_r = [A_0\; A_1\; \cdots \; A_{K-1} ]$.

\subsection{Determining the number of the eigenvalues} \label{S-S-ncs}

 The number of the computed eigenvalues is equal to the number of the leading singular values retained in $\Sigma_0$. Thus a proper truncation in \eqref{eq-svdH} is crucial for obtaining the right eigenpairs. In \cite{AS09, Beyn12} the authors propose to truncate the singular values of $H$ using a predetermined threshold $tol$ such that the first $\bar n_{\mathcal {C}}$ singular values larger than $tol$ are retained in $\Sigma_0$; that is, if $\sigma_1 \ge \cdots \ge \sigma_{\bar n_{\mathcal {C}}} >tol \ge \sigma_{\bar n_{\mathcal {C}}+1} \ge \cdots \ge \sigma_{K\cdot k_S}$, then $V_0 = V(1:\bar n_{\mathcal {C}}, :)$, $\Sigma_0 = \diag( \sigma_1, \cdots, \sigma_{\bar n_{\mathcal {C}}} )$, $W_0 = W(1:\bar n_{\mathcal {C}}, :)$. But in general the optimal threshold $tol$ is affected by several factors, including the conditioning of the matrices in \eqref{eq-ranks}, the distribution of the eigenvalues and the positions of the contour, etc. As a result, it is not easy to choose it beforehand.

Here we suggest two approaches to determine $\bar n_{\mathcal {C}}$ automatically. The first approach is based on the fact that the number of eigenvalues inside $\mathcal {C}$ is equal to the winding number
\begin{equation} \label{eq-ncs}
\bar n_{\mathcal {C}} = {1 \over 2 \pi \im} \int_{\mathcal {C}} \left[{\diff  \over \diff z } \log \det (T_S(z)) \right] \diff z = {1 \over 2 \pi \im} \int_{\mathcal {C}} \mathrm{tr} \left( T_S(z)^{-1} {\diff T_S(z) \over \diff z }  \right)  \diff z,
\end{equation}
where $\mathrm{tr}(\cdot)$ is the trace of matrix. Since in this paper $T_S(z)$ is expressed in the form of \eqref{eq-Zs-expr}, the derivative $\diff T_S(z) / \diff z$ can be computed easily.
The contour integral in \eqref{eq-ncs} can be computed by numerical quadratures, like the one given by \eqref{eq-Ap-num}.

The second approach is based on the gaps in the singular values of $H$ in \eqref{eq-svdH}. We found that the largest gap between two successive singular values denoted by
$$
g_{\max} = \max_{j=1,\cdots, k_S - 1}\left({\sigma_j / \sigma_{j+1}}\right),
$$
typically reaches its maximum at $j=\bar n_{\mathcal {C}}$. This leads to the following strategy:
\begin{equation} \label{eq-gap}
\mathrm{if} \;  g_{\max} \ge \tol_{\gap},  \quad \mathrm{then}\; \bar n_{\mathcal {C}} = \underset{ j=1,\cdots, k_S - 1 }{ \mathrm{Argmax} }  \left({\sigma_j / \sigma_{j+1}} \right),
\end{equation}
where the user-defined constant $\tol_{\gap}$ is introduced only for checking whether a $g_{\max}$ is reasonable; for example, we always set $\tol_{\gap} = 10^{3}$.

In practice, the two approaches can be combined to form more reliable strategies for determining $\bar n_{\mathcal {C}}$. When the moments $A_\alpha$ in \eqref{eq-Ap} are accurately evaluated, \eqref{eq-ncs} and \eqref{eq-gap} should have the same result. Otherwise, if the two approaches result in different numbers of the eigenvalues, one may need to check whether the contour is properly selected; or one can compute the eigenpairs corresponding to both the two numbers and determine which is more reasonable by checking the relative residual $||T_S(\lambda)v||_2/||v||_2$ for each eigenpair $(\lambda, v)$.

\subsection{Algorithm for solving the reduced NEP}
The procedures of the modified block SS method is summarized in Algorithm \ref{alg-TsNEP}. Usually, the contour can be chosen from the ellipses given by the following parameterization
\begin{equation} \label{eq-contour}
\varphi(\alpha) = \gamma + \left [ a \cos(\alpha) + \im  b \sin(\alpha) \right], \quad \alpha \in [0, 2\pi),
\end{equation}
where, $\gamma \in \mathbb{C}$ is the center, $a$ and $b$ are the lengths of the major axis and minor axis, respectively. Then, the moments of $T_S(z)^{-1}$ in \eqref{eq-Ap} can be evaluated by using $N_S$-point trapezoid rule, leading to the following approximation
\begin{equation}\label{eq-Ap-num}
A_\alpha \approx {1 \over  \im N_S} \sum^{N_S-1}_{j=0}  \left[{\varphi(\theta_j)-\gamma \over \rho }\right]^\alpha  \varphi'(\theta_j)  T\left( \varphi(\theta_j) \right)^{-1}, \quad \alpha=0,1,\cdots, 2K-1,
\end{equation}
where $\theta_j = 2\pi  (j+1/2)/N_S, \, j=0, \cdots, N_S$ are the quadrature points.

\begin{algorithm}
\begin{enumerate} [\bf Step 1:]
  \item Choose a contour $\mathcal {C}$, the number of quadrature points $N_S$ on $\mathcal {C}$ and $K$ that determines the number of blocks in Hankel matrices. \label{alg1-step1}

  \item Compute the moments $A_\alpha, \, \alpha=0, \cdots, 2K-1$ from \eqref{eq-Ap}, and the matrices $H$ and $H^<$ \eqref{eq-HH}.

  \item Compute the SVD of $H$ and determine $\bar n_{\mathcal {C}}$ by using the strategy in Section \ref{S-S-ncs}.

  \item Compute the matrix $A$ from \eqref{eq-A}, solve the corresponding standard eigenvalue problem and obtain the $\bar n_{\mathcal {C}}$ eigenpairs $(\lambda, g)$ of the reduced NEP from \eqref{eq-eigens}.
  \vspace{-8pt}
\end{enumerate}
\caption{Modified block SS algorithm for solving the reduced NEP \eqref{eq-Ts-NEP}. \label{alg-TsNEP}}
\end{algorithm}

We remark that the modified block SS method is also based on Theorem \ref{thm1} in Section \ref{S-basics}. In particular, following the block SS method it is possible to solve the original NEP \eqref{gnep} using the moments $\hat M_\alpha$ computed in \eqref{eq-malpha-qud}. However, as we have pointed out in the introduction, this approach can be unstable and inaccurate when $L$ is not very large. Bear in mind that a large $L$ implies a high computational burden in large scales. But the situation is quite different in small scales. Since the computational cost is negligible, one can afford to choose $L$ to be the dimension of the matrix $T_S$, and at the same time to use high order quadratures for the contour integration. As a result, the modified block SS method is very stable and accurate and thus an efficient approach for solving NEPs of small and moderate sizes.

\section{Summary of the proposed \cirrs{} algorithm} \label{S-algorithm}

The \cirrs{} method proposed in this paper can be implemented following Algorithm \ref{alg-cirrip}.

Step \ref{alg-cirrip-step1} is the initialization phase. The foremost task is to choose a suitable $\mathcal {C}$ within which the eigenvalues are searched for. For many engineering problems, one may have a-priori information about the locations of the interested eigenvalues from experiments, theoretical predictions, eigenvalues of similar structures, etc.
Another cheaper way to acquire such information is to perform an analysis on a coarse mesh or a simplified model, both of which have been frequently used in literature; see e.g., \cite{Voss07}.
The contour is usually chosen to be an \emph{ellipse} \eqref{eq-contour} or a \emph{rectangle} enclosing the interested eigenvalues. The sampling points can then be set as the quadrature points on the contour.

Given a value of $L$, the number $N$ should be chosen such that the ratio of the first to the last singular values of $\hat S$ is large enough, e.g., larger than $10^{14}$. This is often the case when $N \cdot L$ is 2 or 3 times larger than the number of the eigenvalues $\bar n_{\mathcal {C}}$; see \eqref{eq-rank-nc}. The number $L$ has to be at least not less than the algebraic multiplicity of the interested eigenvalues. In some cases (e.g., structural modal analysis) this can be estimated from the symmetry of the problems. Otherwise, one can use a large $L$, but this would increase the computational cost to a certain extent, depending on the properties of the matrices.

In Steps \ref{alg-cirrip-step2} and \ref{alg-cirrip-tsvd} we construct the eigenspaces. The method has been detailed in Section \ref{S-eigenspace-samp}. The threshold $\delta$ of the truncated SVD can always be set as $\delta = 10^{-14}$. In Steps  \ref{alg-cirrip-step4} and \ref{alg-cirrip-step5} we solve the reduced NEP and compute the eigenpairs of the orginal NEP using the method in Section \ref{S-reducedNEP}. Note that the same contour is used in solving the reduced NEP and constructing the eigenspaces.

Finally, the \cirrs{} method can be easily implemented in conjunction with other software. In fact, only Steps \ref{alg-cirrip-step2} and \ref{alg-cirrip-step4} involve the operations with the matrix $T$ and can be carried out by the host software. The other steps can be considered as preprocessing and postprocessing parts for the host software. In this paper, the \cirrs{} is implemented into the FEM software \textsc{Ansys}\textregistered{} and our fast BEM program; see the next section for numerical results. The \cirrm{} algorithm can be implemented analogously as Algorithm \ref{alg-cirrip}. In this situation, the order of moments $K'$ has to be initialized in Step \ref{alg-cirrip-step1}, and the matrix $\hat S$ in Step \ref{alg-cirrip-tsvd} has to be replaced by $\hat M$.

\begin{algorithm}
\begin{enumerate} [\bf Step 1:]
  \item Choose a contour $\mathcal {C}$ and the parameters $N$ and $L$, choose a random matrix $U \in \mathbb{C}^{n \times L}$, and determine the sampling points $z_i, \, i=0, \cdots, N-1$ on $\mathcal {C}$. \label{alg-cirrip-step1}

  \item Compute $T(z_i)^{-1} U, \, i=0, \cdots, N-1$. \label{alg-cirrip-step2}

  \item Formulate $\hat S$ by \eqref{eq-hatS}. Compute the truncated SVD $\hat S \approx S \Sigma_0 W^H_0$ such that the first $k_S$ singular values larger than $\delta\cdot\sigma_1$ are retained. \label{alg-cirrip-tsvd}

  \item Compute $T_S(z) = S^H T(z)S$, and solve the reduced NEP \eqref{eq-Ts-NEP} by Algorithm \ref{alg-TsNEP} to obtain $\bar n_{\mathcal {C}}$ eigenpairs $(\lambda_j, g_j), \, j=1, \cdots, \bar n_{\mathcal {C}}$. \label{alg-cirrip-step4}

  \item Compute the eigenpairs $(\lambda_j, v_j)$ of the original NEP \eqref{gnep} via $(\lambda_j, S g_j), \, j=1, \cdots, \bar n_{\mathcal {C}}$.  \label{alg-cirrip-step5}
  \vspace{-8pt}
\end{enumerate}
\caption{The \cirrs{} algorithm. \label{alg-cirrip}}
\end{algorithm}

\section{An extension to NEPs in the BEM} \label{S-bem}

The \cirrs{} algorithm can be readily applied when the matrix $T(z)$ is explicitly expressed as
\begin{equation} \label{eq-Z-expr}
T(z)  = \sum_{j=1}^J T_j f_j(z),
\end{equation}
where $f_j(z)$ and $T_j$ are certain functions and the coefficient matrices, because in this case the expression of the reduced matrix $T_S(z)$ can be easily obtained as
\begin{equation} \label{eq-Zs-expr}
T_S(z)  = \sum_{j=1}^J  \left(S^H T_j S \right)f_j(z).
\end{equation}
By assuming that the number of the terms $J$ is much less than $n$ (for example, less that 100), the coefficient matrices $S^H T_j S$ can then be computed with an affordable computational cost. Note that the evaluation of $S^H T_j S$ is suitable for parallelization.

However, there is a more general case where the matrix $T(z)$ is not inherently of the form \eqref{eq-Z-expr}. A typical example is the NEPs in the BEM, in which the entries $T_{ij}(z)$ of $T(z)$ are different functions of $z$, and the explicit expressions of $T_{ij}(z)$ are difficult to obtain when singular integrations are involved. Consequently, the \cirrs{} can not be directly applied.
In this section, the NEPs of the BEM in acoustics is taken as a model problem. After a concise review of the basic formulation of the NEPs, the Chebyshev interpolation of $T(z)$ is used to convert it into the form \eqref{eq-Z-expr}.

\subsection{NEPs of BEM in acoustics} \label{S-S-bemexs}

The steady-state acoustic wave pressure $p$ at any location $\bbx \in R^{3}$ in a bounded domain $D$, enclosed by a boundary surface $\partial D$, is governed by the second-order Helmhotz partial differential equation
\begin{equation}\label{hpde}
\nabla^2 p  + z^2 p  = 0,
\end{equation}
where, $z$ is used to denote the wavenumber, $z = \omega / c = 2 \pi f /c$, and $c$ is the speed of sound. In order to uniquely define the pressure field in the domain $D$, proper boundary conditions, e.g., the Dirichlet or Neumann boundary conditions, have to be specified on the closed boundary surface $\partial D$.

The BEM does not treat \eqref{hpde} directly, instead, it solves the following Helmholtz boundary integral equation (BIE)
\begin{equation}\label{hbie}
C(\bbx) p(\bbx)  + \int_{\partial D}{ \partial F(\bbx, \bby) \over \partial n_{\bby} } p(\bby) \diff \bby = \int_{\partial D} F(\bbx, \bby) { \partial p(\bby) \over \partial n_{\bby} }  \diff \bby,
\end{equation}
which can be obtained by using the Green's second identity to the Helmholtz PDE. In \eqref{hbie}, $C(\bbx)$ is the free term coefficient at $\bbx$, $n_{\bby}$ denotes the outward normal at $\bby$, $F(\bbx, \bby) = \exp(\im z |\bbx-\bby|) / |\bbx-\bby|$ is the fundamental solution of the Helmholtz equation in three dimensions.

In this paper, the Nystr\"om method in \cite{CWX15} is used to discretize the Helmhotz BIE \eqref{hbie}. The boundary surface $\partial D$ is discretized into curved triangular quadratic elements. The 6-node conforming quadratic element is used in the geometry approximation, and the 6-node non-conforming quadratic element is used for the approximation of the unknown functions; the nodes are the projections of the Gauss quadrature points on the reference element. The Nystr\"om discretization leads to the following standard system in the BEM,
$$
H(z)u(z) = G(z)q(z),
$$
where $H$ and $G$ are complex square matrices obtained by the Nystr\"om discretization of the left and right side of BIE \eqref{hbie}, respectively; $u$ and $q$ are vector collections of the nodal pressure $p$ and its normal direvative $\partial p / \partial n$ which is related to the normal velocity of the acoustic media.

By considering the boundary conditions and gathering together the unknown quantities into a vector $v$, one obtains a linear system of form $T(z) v(z) = b(z)$, where the vector $b$ accounts for the contribution of the boundary conditions and should be set to zero in eigenvalue problems. Therefore, the NEPs in the BEM come into the general form \eqref{gnep},
$$T(z) v  = 0.$$
In accordance with \eqref{gnep}, the dimensions of all the matrices and vectors mentioned here are denoted by $n$.

As the Nystr\"om discretization based on non-conforming elements is used, the entries of $H$ and $G$ associated with well-separated nodes $\bbx_i$ and $\bby_j$ are just the values of double- and single-layer integral kernels, that is,
\begin{equation} \label{eq-nyes}
H_{ij}(z) = {\partial \over \partial n_{\bby_j}} {\exp(\im z r_{ij}) \over r_{ij}}, \quad G_{ij}(z) =  {\exp(\im z r_{ij}) \over r_{ij}},
\end{equation}
where $r_{ij} = |\bbx_i-\bby_j|$. When nodes $\bbx_i$ and $\bby_j$ are not well-separated, i.e., the element integral corresponds to $\bbx_i$ and the boundary element of $\bby_j$ is singular or nearly singular, the expressions of $H_{ij}(z)$ and $G_{ij}(z)$ in terms of $z$ could be much more complicated. As a result, the entries of $T(z)$ are generally distinct functions of $z$, and for some of them the expressions are very difficult to obtain. This is an intrinsic feature of the BEM. The situations in the collocation and Galerkin discretizations are even more troublesome.

In order to enable large-scale and wide frequency band solutions, the BEM in this paper is accelerated by using the fast directional multipole method; see \cite{CWX15}.

\subsection{Chebyshev interpolation of BEM matrix} \label{S-S-chebyshev}

Here by using Chebyshev polynomial interpolation, the matrix $T$ is transformed into the form \eqref{eq-Z-expr} so that the reduced matrix $T_S(z)$ can be computed with moderate computational cost.

To simplify the presentation, the domain of interest is assumed to be the interval $[-1,\,1]$. Let
\begin{equation} \label{eq-T-interp}
T(z) \approx P(z) = P_0 \tau_0(z) + \cdots + P_d \tau_d(z)
\end{equation}
be a polynomial interpolation of the matrix $T(z)$ of degree $d$ that satisfies
$$
T(x_j) = P(x_j), \quad j=0, \dots, d,
$$
where $\tau_j$ are Chebyshev polynomials of order $j$, $x_j \in [-1,\,1]$ are the Chebyshev points of the first kind,
$$
x_j = \cos (\theta_j), \quad \theta_j =  {j+ {1\over 2} \over d+1} \pi,
$$
 and $P_j$ are coefficient matrices that can be computed by
\begin{equation} \label{eq-T-coefs}
 P_j = {2 \over d+1} \sum_{k=0}^{d} T(x_k) \cos( {j\theta_k}).
\end{equation}

 By combining \eqref{eq-T-interp}, \eqref{eq-T-coefs} and \eqref{eq-Ts-NEP}, one obtains
 \begin{equation} \label{eq-Ts-interp}
T_S(z) \approx  P_S(z) := S^H P(z) S = \hat P_0 \tau_0(z) + \cdots + \hat P_d \tau_d(z),
\end{equation}
where the coefficient matrices $\hat P_j \in \mathbb{C}^{k_S \times k_S}$ are given by,
\begin{equation} \label{eq-Ts-coefs}
\hat P_j = {2 \over d+1} \sum_{k=0}^{d} \left[S^H T(x_k) S\right] \cos( {j\theta_k}).
\end{equation}
Thus, the coefficient matrices $\hat P_j$ of the Chebyshev interpolation $P_S(z)$ can be computed by the matrix-vector multiplications with the original system matrix $T$ at Chebyshev points $x_k$. These operations can be efficiently parallelized.

As $T$ is a holomorphic function on a neighborhood of the real interval $[-1,\,1]$, the error of the Chebyshev interpolant in \eqref{eq-T-interp} decays exponentially with the degree $d$ \cite{EK12}. Thus, a moderate $d$ is often sufficient to ensure good accuracy.

\section{Numerical examples} \label{S-ne}

The performance of the \cirrs{} algorithm is demonstrated by several representative examples, including some standard benchmark problems for testing new nonlinear eigensolvers in literature, and some practical problems in engineering finite element and boundary element analyses. Comparisons with the state-of-the-art nonlinear eigensolvers, such as the \cirrm{} and the NLEIGS in \cite{guttel2014nleigs} are also conducted.

All the computations were performed on a personal computer with an Intel$^\circledR$
Core$^\mathrm{TM}$ i3-2100 (3.10 GHz) CPU and 16 GB RAM, except the example in Section \ref{S-S-paf} which was carried out on a Server with eight 8-core Intel Xeon E7-8837 (2.67GHz) processors and 256 GB RAM.

\subsection{Three benchmark examples}\label{S-ne-4bench}

First, we compare the accuracy of the \cirrs{} and \cirrm{} by using three benchmark examples from the NLEVP collection \cite{BHM13}.
The accuracy of an eigenpair $(\lambda, v)$ is measured by the relative residual $||T(\lambda)v||_2/||v||_2$.
The first two examples belong to the quadratic eigenvalue problem which is also nonlinear. The coding and computation are conducted using Matlab R2009a. The linear systems in computing $T(z_i)^{-1} U$ are solved using the Matlab backslash operator \verb"`\'" in sparse mode.

\subsubsection{Description of the examples}

\begin{enumerate}[(1)]

 \item {\it Acoustic-1D.} This is a quadratic eigenvalue problem of the form $T(\lambda)v = (\lambda^2 M + \lambda C + K_{\rm{s}} ) v = 0$, arising from the finite element discretization of the time harmonic
wave equation \eqref{hpde} in a one-dimensional interval $[0, 1]$, where the boundary conditions are Dirichlet ($p=0$) at one end and impedance ($q + {2\pi f \im \over \zeta}p = 0$) at another end. The $n \times n$ matrices are defined by
$$
     M = -{4\pi^2 \over n}\left(I_{n} - {1 \over 2} e_n e_n^T \right), \quad C = {2\pi f\im \over \zeta} e_n e_n^T,
     \quad K_{\rm{s}} = {n}\begin{bmatrix}
  2 & -1 &   &   \\
  -1 & \ddots & \ddots &   \\
    & \ddots & 2 & -1 \\
    &   & -1 & 1 \\
\end{bmatrix}.
     $$
 We use command \verb"nlevp(`acoustic_wave_1d', 1000, 1)" to generate these matrices with size $n = 1000$, impedance $\zeta=1$ and wave speed $c=1$.

 \item {\it Concrete.} This quadratic eigenvalue problem is of the form $T(\lambda)v = (\lambda^2 M + \lambda C + (1+\im \mu)K_{\rm{s}} ) v = 0$, arising in a model of a concrete
structure supporting a machine assembly \cite{FPS00}. The matrix $M$ is real diagonal and of low rank. The viscous damping matrix $C$ is pure imaginary and diagonal. The matrix $K_{\rm{s}}$ is complex symmetric, and the factor $(1+\im \mu)$ adds uniform hysteretic damping. These matrices are generated by \verb"nlevp(`concrete', 0.04)" with size $n = 2472$ and $\mu = 0.04$.

  \item {\it String.} This is a rational eigenvalue problem arising in the finite element discretization of a boundary eigenvalue problem
describing the eigenvibrations of a string with a load of unit mass attached by an elastic spring of unit stiffness \cite{SSI06}. It has the form
$T(\lambda)v =  (T_1 + {\lambda \over \lambda-1} e_n e_n^T - \lambda T_3 ) v = 0$, where, $n$ is the number of equally-spaced elements and
$$
T_1 = {n}\begin{bmatrix}
  2 & -1 &   &   \\
  -1 & \ddots & \ddots &   \\
    & \ddots & 2 & -1 \\
    &   & -1 & 2 \\
\end{bmatrix},
\quad
T_3 = {1 \over 6n}\begin{bmatrix}
  4 & 1 &   &   \\
  1 & \ddots & \ddots &   \\
    & \ddots & 4 & 1 \\
    &   & 1 & 2 \\
\end{bmatrix}.
$$

We use \verb"nlevp(`loaded_string', 5000, 1, 1)" to generate the matrices and compute the $32$ eigenvalues in the interval $[3, 10000]$. The size of this problem is $n=5000$.
\end{enumerate}

The controlling parameters of the \cirrs{} and \cirrm{} used in the three examples are listed in Table \ref{tab-parameters}.
The elliptical contours \eqref{eq-contour} are used, and the $N$-point trapezoidal quadrature rule \eqref{eq-Ap-num} is employed to determine the sampling points and compute the moments.
The orthogonal basis of the eigenspace is computed by using SVD. We do not truncate the SVD so that the dimensions of the resultant eigenspaces corresponding to $\hat S$ and $\hat M$ are maximal and have the same dimension. This should reflect the full power of the two eigenspaces respectively. To ensure the rank condition \eqref{eq-ranks} and high accuracy in solving the reduced NEPs, proper values of $K$ and $N_S$ are used; see the last two columns of Table \ref{tab-parameters}.

\begin{figure}[htb]
\centering
  \epsfig{figure=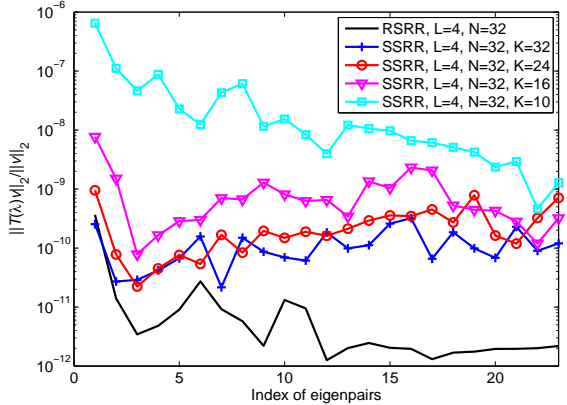,width=0.5\textwidth}
\caption{Influence of $K'$ on the accuracy of \cirrm{} and the comparison with \cirrs{} using the {\it `String'} example.}
\label{fig-string-NLK}
\end{figure}

Note that in the comparisons with the \cirrm{} algorithm, we always set $K'=N$, because this will result in the largest subspaces for the moment scheme. For a given $L$, using a smaller $K'$ will lead to a smaller subspace and thus inferior eigensolutions; from the viewpoint of the computational efficiency, this implies that we have not made full use of the probed resolvent $T(z)^{-1}U$ at the sampling points. This statement can be verified by Figure \ref{fig-string-NLK}, which illustrates the relative residuals of the computed eigenpairs from an experiment using the {\it `String'} example. We have used the contour in Table \ref{tab-parameters}, set $L=4$ and $N=32$, and increased $K'$ of \cirrm{} from 10 to 32. When showing the residuals, the eigenvalues are first sorted in ascending order according to their absolute values, then the residual of each eigenpair is shown against its index in the sorted list.

\begin{table}[!h]
\begin{center}
\caption{Values of the controlling parameters of the \cirrs{} and \cirrm{} algorithms in the three benchmark examples}\label{tab-parameters}
\vspace{-.8\baselineskip}
\begin{tabular*}{\textwidth}{@{\extracolsep{\fill}}cccccccc@{}}\toprule
%&\multicolumn{3}{c}{rre}    &\multicolumn{2}{c}{Present}\\
                   &  $\gamma$    & $(a,b)$      & $L$  & $N$    & $K'$    & $N_S$ & $K$\\
\hline
\emph{Acoustic-1D} & $9.9+0.8\im$ & $(10.1,1.01)$& 2    & 100    & 100     & 1000             & 2\\
\emph{Concrete}    & $-1+52.5\im$ & $(5.25,52.5)$& 2    & 100    & 100     & 1000             & 2\\
\emph{String}      & $5001.5$     & $(4998.5,249.925)$& 1    & 100    & 100     & 1000         & 8\\
\bottomrule
\end{tabular*}
\end{center}
\end{table}

\subsubsection{Numerical results and discussions}

\begin{figure}
\centering
\begin{subfigure}[t]{.5\textwidth}
  \centering
  \epsfig{figure=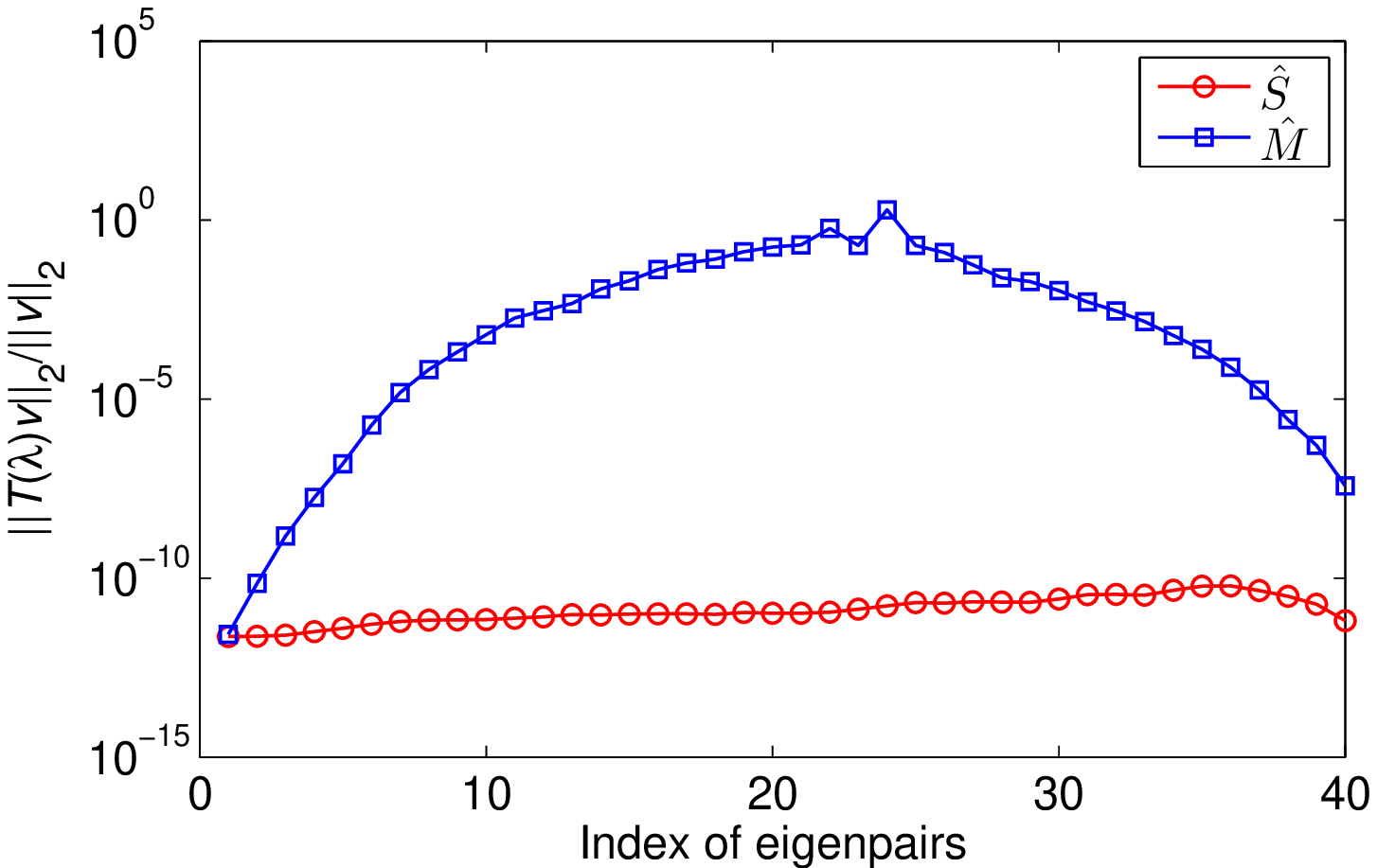,width=0.9\textwidth}
  \caption{Residuals of eigenpairs, \emph{Acoustic-1D} example}
  \label{fig-acoustic-err}
\end{subfigure}%
\begin{subfigure}[t]{.5\textwidth}
  \centering
  \epsfig{figure=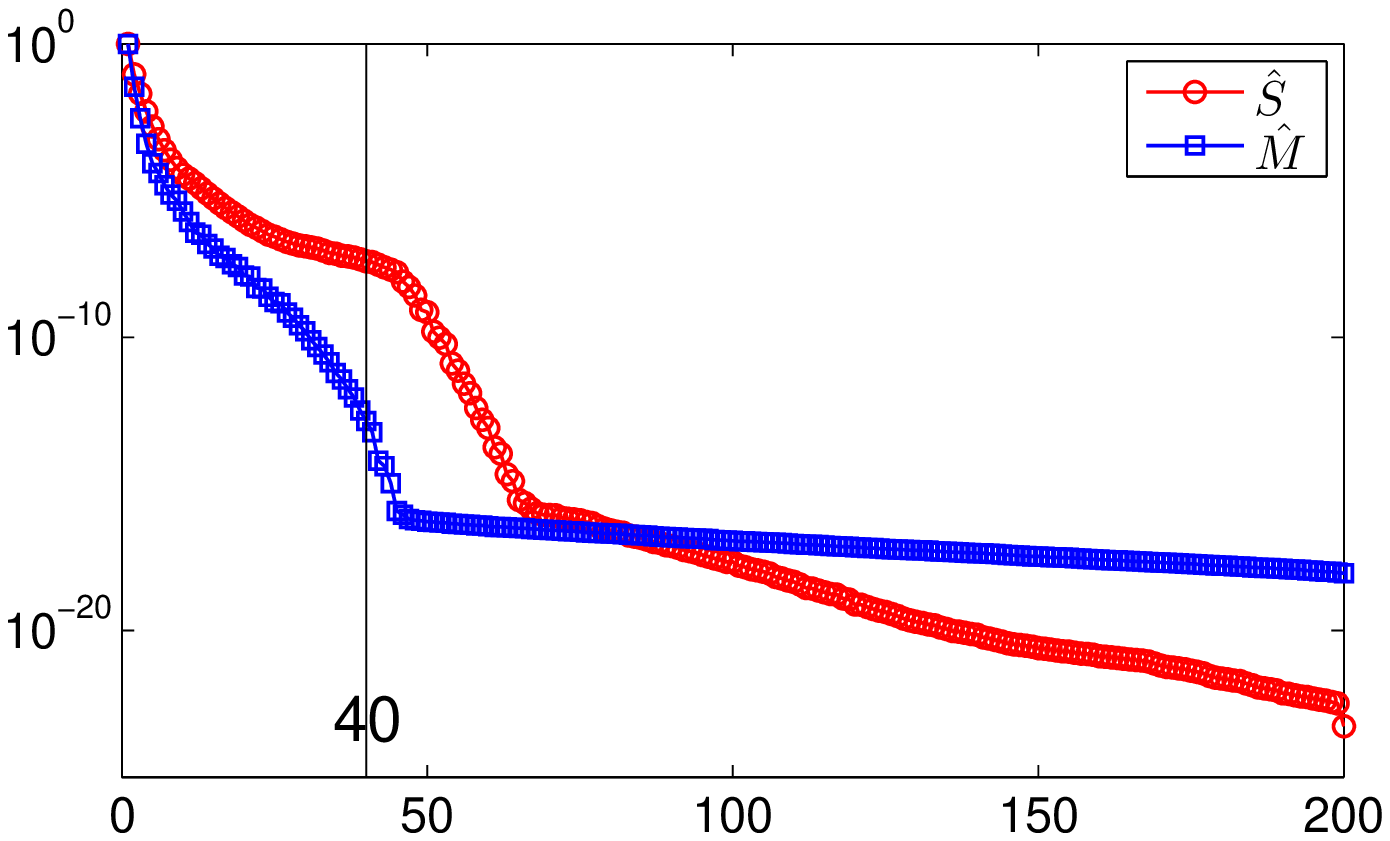,width=0.9\textwidth}
  \caption{Singular values of $\hat S$ and $\hat M$, \emph{Acoustic-1D} example}
  \label{fig-acoustic-sv}
\end{subfigure}\\

\begin{subfigure}[t]{.5\textwidth}
  \centering
  \epsfig{figure=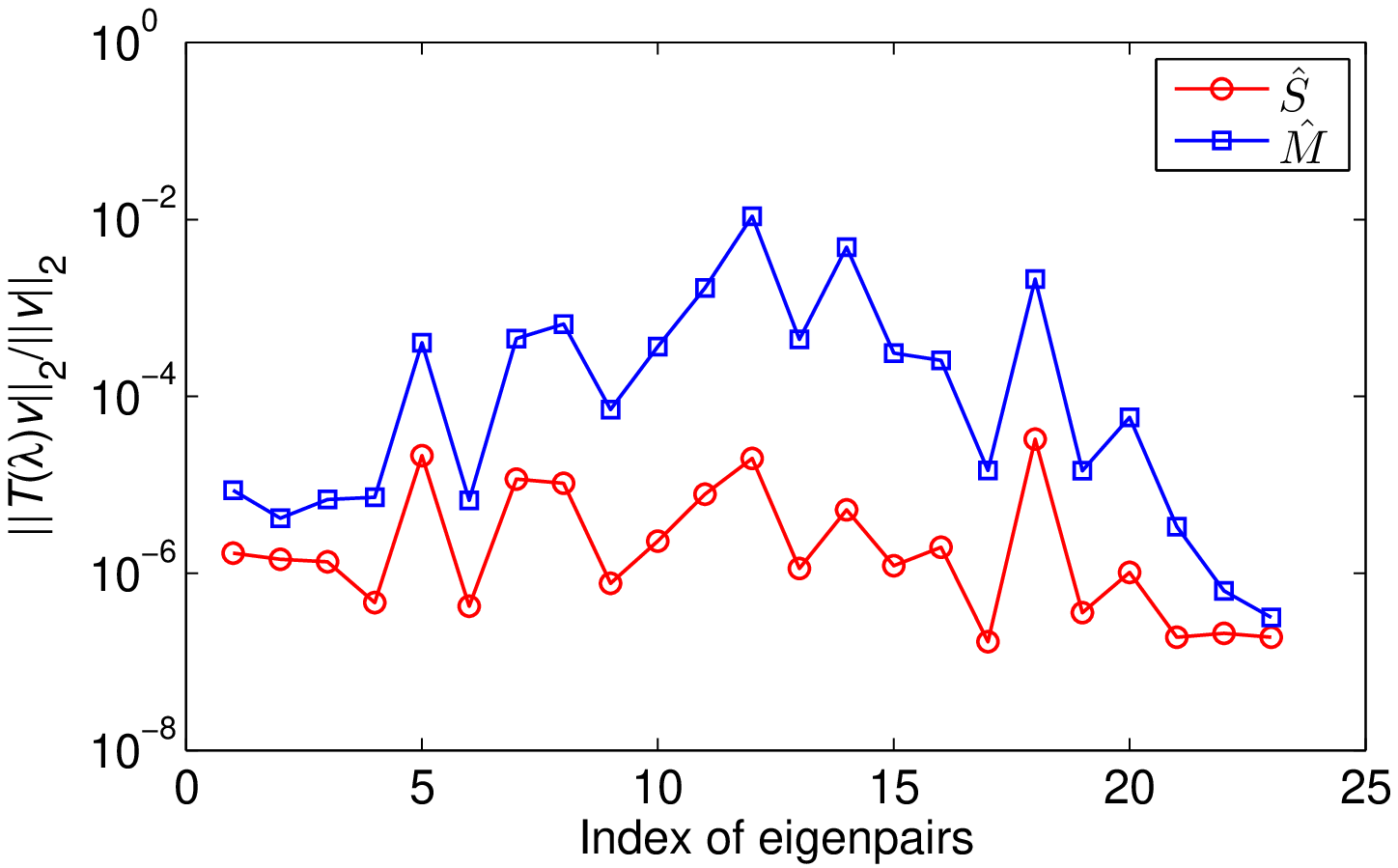,width=0.9\textwidth}
  \caption{Residuals of eigenpairs, \emph{Concrete} example}
  \label{fig-conc-err}
\end{subfigure}%
\begin{subfigure}[t]{.5\textwidth}
  \centering
  \epsfig{figure=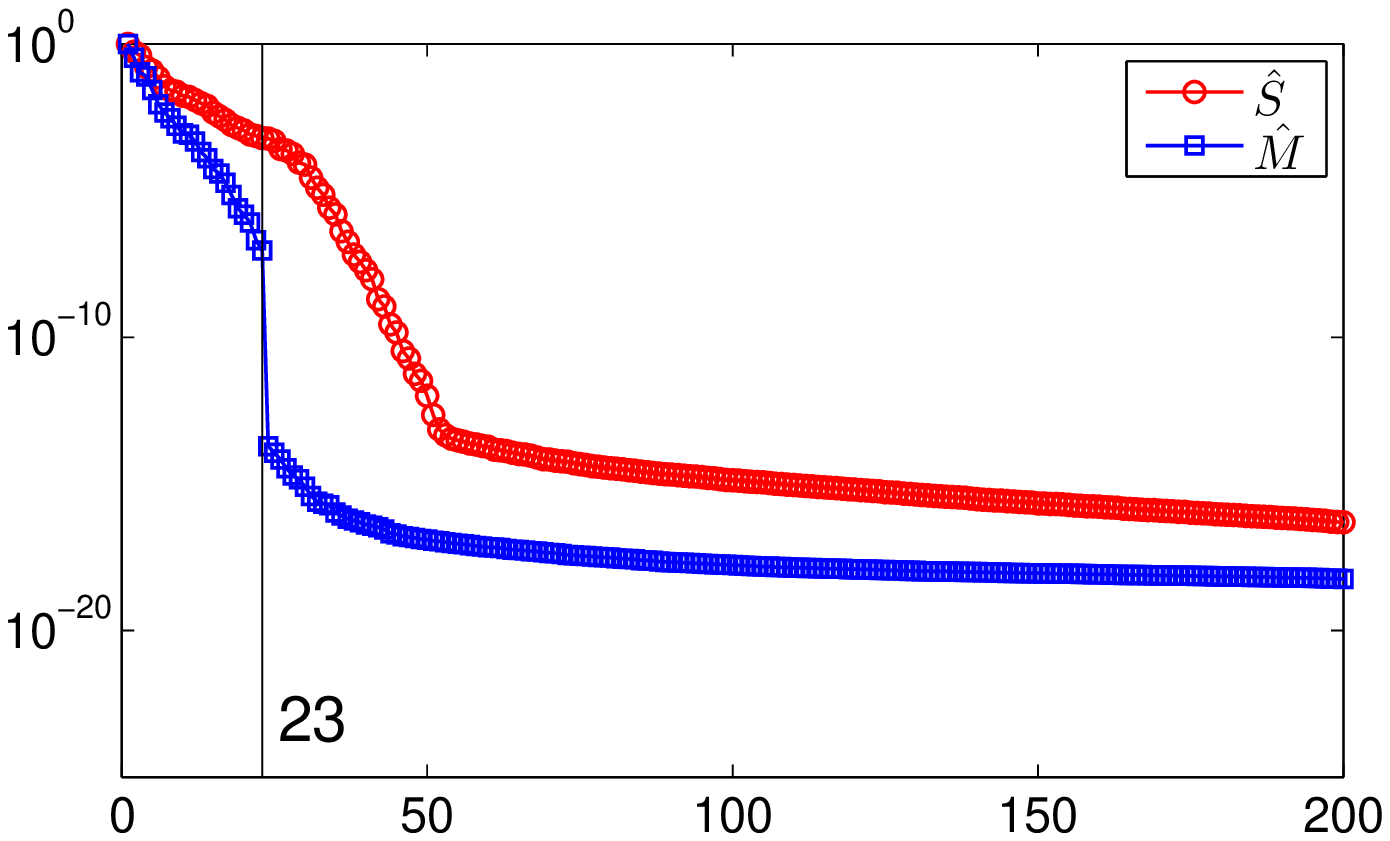,width=0.9\textwidth}
  \caption{Singular values of $\hat S$ and $\hat M$, \emph{Concrete} example}
  \label{fig-conc-sv}
\end{subfigure}\\

\begin{subfigure}[t]{.5\textwidth}
  \centering
  \epsfig{figure=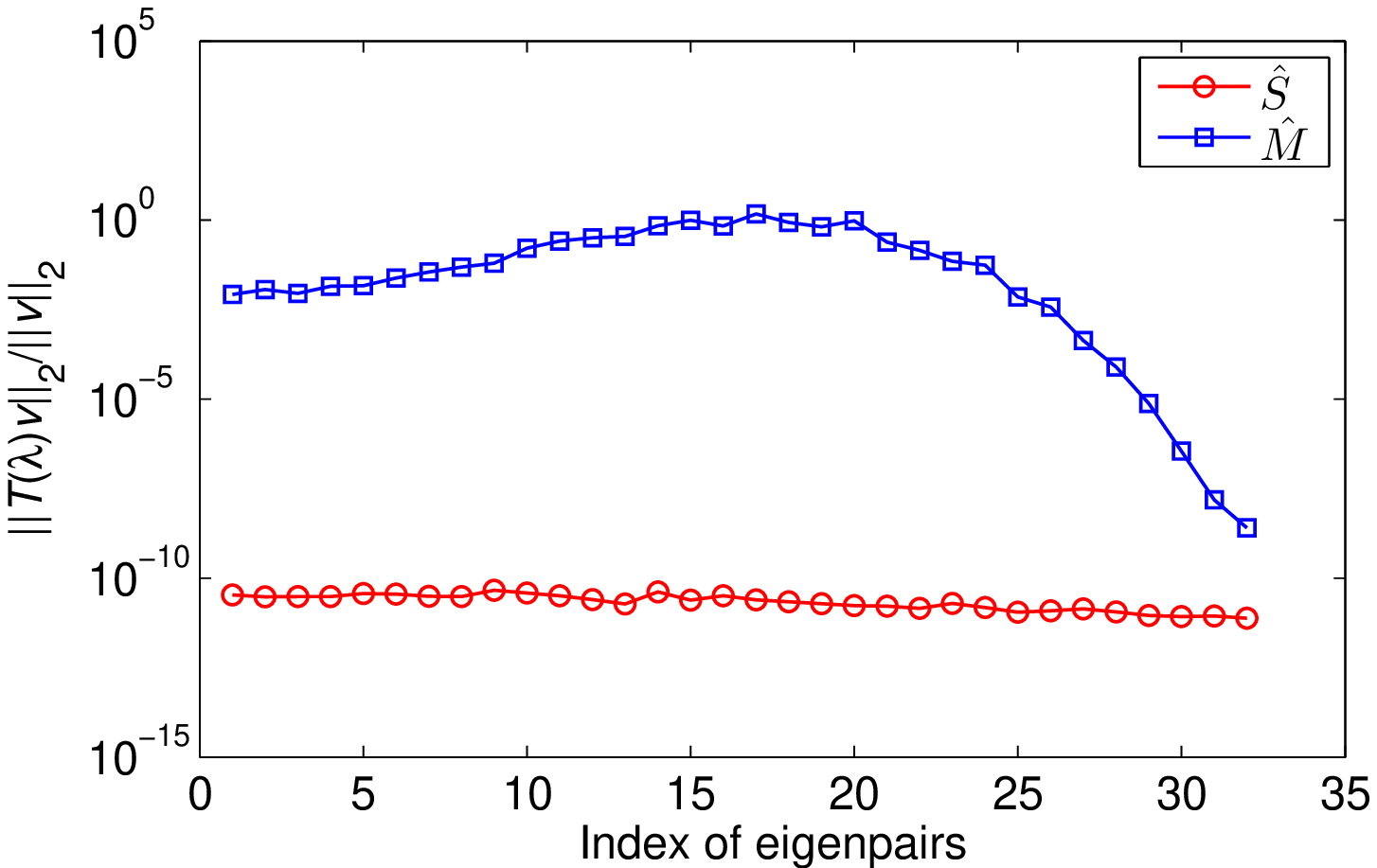,width=0.9\textwidth}
  \caption{Residuals of eigenpairs, \emph{String} example}
  \label{fig-string-err}
\end{subfigure}%
\begin{subfigure}[t]{.5\textwidth}
  \centering
  \epsfig{figure=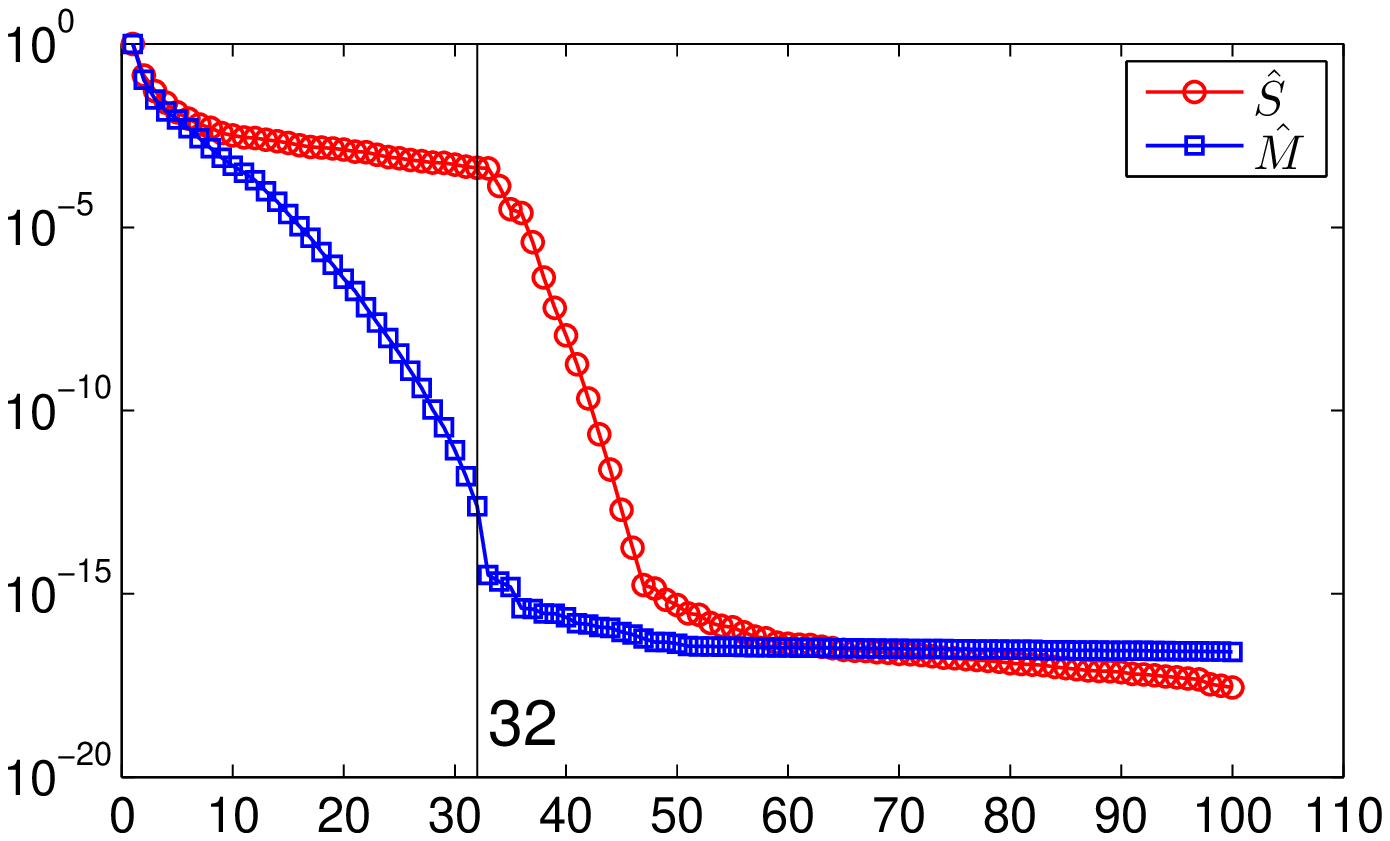,width=0.9\textwidth}
  \caption{Singular values of $\hat S$ and $\hat M$, \emph{String} example}
  \label{fig-string-sv}
\end{subfigure}
\caption{Comparisons of the sampling scheme and the moment scheme for constructing eigenspaces for the three benchmark examples.}
\label{fig-four-benchs}
\end{figure}

\begin{figure}
\centering
\begin{subfigure}[t]{.5\textwidth}
  \centering
  \epsfig{figure=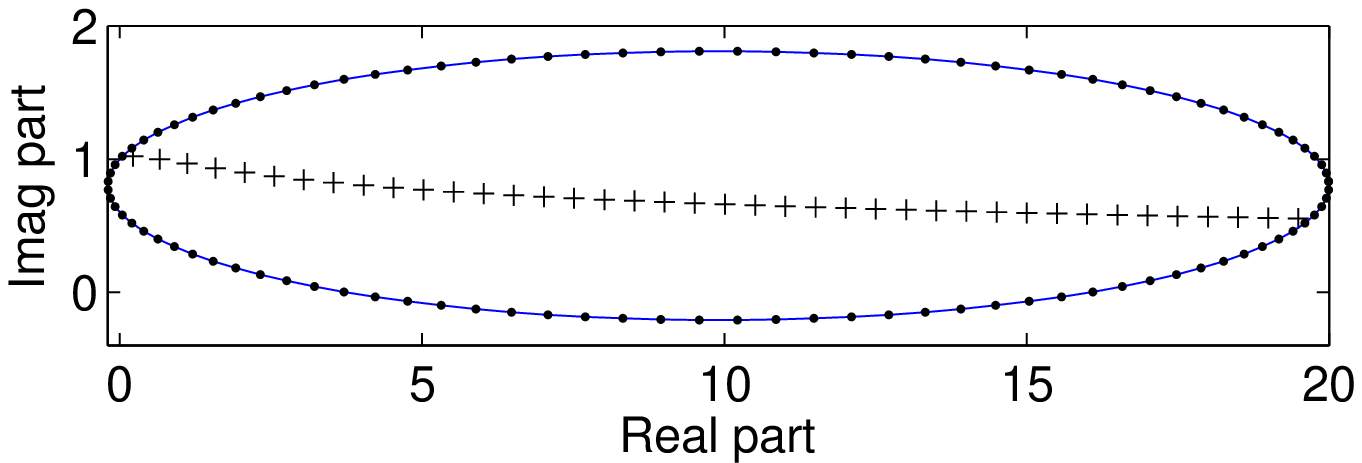,width=0.975\textwidth}
  \caption{\emph{Acoustic-1D} example}
  \label{fig-acoustic-eigs}
\end{subfigure}\\

\begin{subfigure}[t]{.5\textwidth}
  \centering
  \epsfig{figure=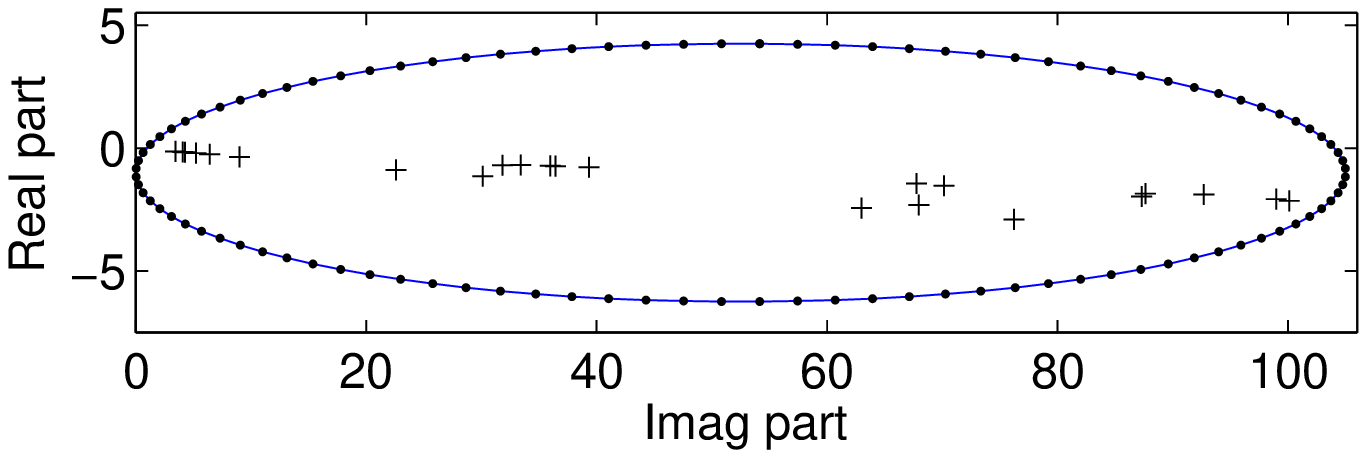,width=0.975\textwidth}
  \caption{\emph{Concrete} example}
  \label{fig-conc-eigs}
\end{subfigure}%
\begin{subfigure}[t]{.5\textwidth}
  \centering
  \epsfig{figure=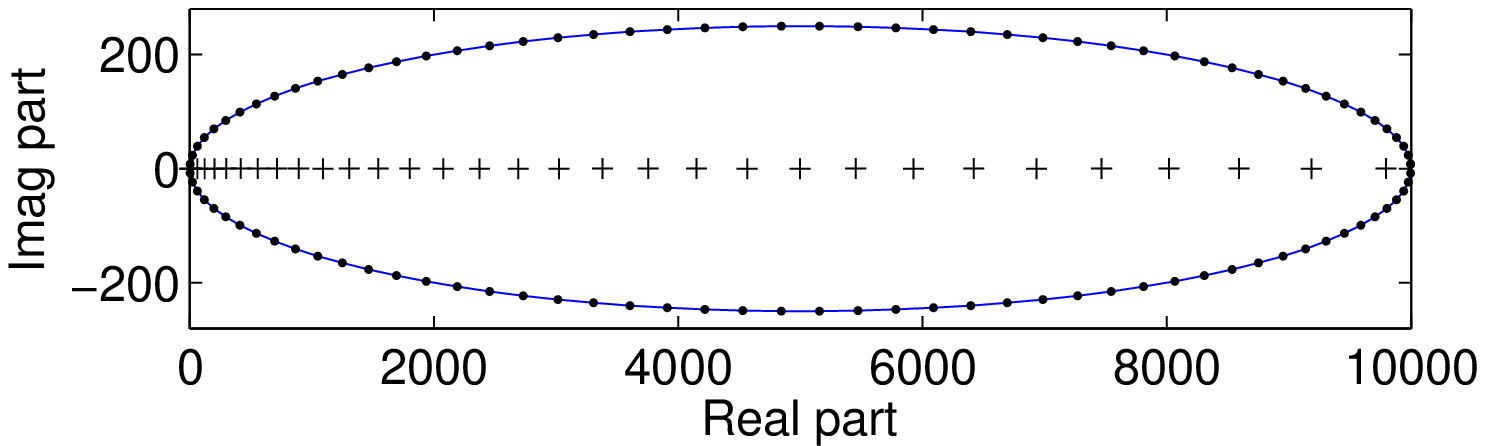,width=0.975\textwidth}
  \caption{\emph{String} example}
  \label{fig-string-eigs}
\end{subfigure}
\caption{Computed eigenvalues ($\scriptstyle\pmb{+}$) by using $\hat S$, contours ($-$) and sampling points ($\bullet $) of the three benchmark examples. Note that in Figures \ref{fig-conc-eigs} $x$-axis corresponds to the imaginary part while $y$-axis corresponds to the real part, which are different from Figures \ref{fig-acoustic-eigs} and \ref{fig-string-eigs}.}
\label{fig-four-benchs-eigs}
\end{figure}

The numerical results are summarized and illustrated in Figure \ref{fig-four-benchs}. For each example, the relative residuals of the computed eigenpairs and the singular values of the matrices $\hat S$ and $\hat M$ are shown in the same row. Clearly, the residuals of eigenpairs computed by the \cirrs{} are typically several orders lower than the residuals obtained by the \cirrm{}.
As explained in Section \ref{S-eigenspace-mom}, this remarkable difference is caused by the distinct quality of the eigenspaces. To be specific, consider the behavior of singular values of $\hat S$ and $\hat M$, as shown in plots (\ref{fig-acoustic-sv}), (\ref{fig-conc-sv}) and (\ref{fig-string-sv}) of Figure \ref{fig-four-benchs}, where the singular values are scaled so that the first one is $1$. We use the vertical line in each plot to indicate the number of the eigenvalues $\bar n_\mathcal {C}$. We see that, typically, the singular values of $\hat M$ decrease much faster than those of $\hat S$. In plots (\ref{fig-acoustic-sv}) and (\ref{fig-string-sv}), the singular values of $\hat M$ drop to around $10^{-15}$ at the $\bar n_\mathcal {C}$-th singular value, indicating that the rank condition \eqref{eq-rank-M-KL} cannot be satisfied with good accuracy.  Whereas, the situation for $\hat S$ is much better, and an obvious decay is observed only after the $\bar n_\mathcal {C}$-th singular value. An exception can be seen from plot (\ref{fig-conc-sv}), where the singular values of $\hat M$ do not decrease as fast as in the other two examples, but still faster than $\hat S$. This may explain why the discrepancy between the residuals of $\hat M$ and $\hat S$ in this example is not so remarkable as in the other two examples.

In summary, the numerical results confirm the basic claim in Section \ref{S-eigenspace-mom}: the use of the high order moments in $\hat M$ may lead to rank-deficiency and finally fail to ensure the rank condition \eqref{eq-rank-M-KL} for proper extraction of all the eigenvalues. On the contrary, the matrix $\hat S$ can always satisfy the condition \eqref{eq-rank-nc} securely, leading to high quality eigenspaces. Our computation also indicates that, to achieve good eigenspaces using $\hat S$, the sampling points should be close to the eigenvalues so that they can be effectively probed. This is why the contours in these examples are chosen to be oblate ellipses; see the column of $(a,b)$ in Table \ref{tab-parameters}. Figure \ref{fig-four-benchs-eigs} exhibits the locations of the sampling points and the computed eigenvalues by the \cirrs{}.

Finally, we notice that in all the three examples the values of $L$ and $N$ in Table \ref{tab-parameters} are chosen to ensure that the correct number of the eigenvalues $\bar n_\mathcal {C}$ can be obtained by using the \cirrm{}. As we have seen, in this situation the \cirrs{} always performs much better than the \cirrm{}. But those values are by no means optimal for the \cirrs{}. For example, we have used $L=2$ in the \emph{Acoustic-1D} example before. If we use $L=1$, the 40 correct eigenpairs can still be accurately extracted by using the \cirrs{} with almost the same accuracy as shown in Figure \ref{fig-acoustic-err}. However, when using the \cirrm{}, only 32 ``eigenpairs'' can be obtained and most of them are incorrect.

\subsection{Comparisons with other methods}\label{S-ne-comparison}

Here we compare the performance of the \cirrs{} with the recently developed NLEIGS method in \cite{guttel2014nleigs}, which has been shown to be a promising and robust nonlinear eigensolver based on linearization. Our comparison uses the \textit{`gun'} problem of the NLEVP collection, which has been extensively used to test the performance of nonlinear eigensolvers. In this problem the matrix $T(z)$ is of the form
$$
T(z) = K_{\rm{s}}- z^2 M + \im \sqrt{z^2-\sigma_1^2} \,W_1 + \im \sqrt{z^2-\sigma_2^2} \,W_2,
$$
with $K_{\rm{s}}$, $M$, $W_1$ and $W_2$ being real symmetric matrices of the size $9956 \times 9956$. We take $\sigma_1 = 0$ and $\sigma_2 = 108.8774$.

This example has already been included in the NLEIGS code \verb"`nleigs_v0.5'" (available from
\url{http://twr.cs.kuleuven.be/research/software/nleps/nleigs.php}), where the 21 eigenvalues in the upper half-disk with center 223.6 and radius 111.8 are computed using four different variants of the NLEIGS. Here we compare with the two most efficient variants: \emph{Variant R2} and \emph{Variant S}.

\begin{figure}[htb]
\centering
  \epsfig{figure=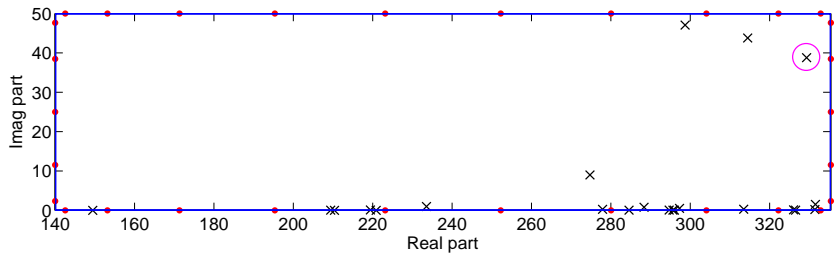,width=0.8\textwidth}
\caption{The eigenvalues ($\scriptstyle\pmb{\times}$), contour (solid line) and the sampling points ({\color{red}$\bullet $}, $N=30$) of the gun example.}
\label{fig-gun-ecs}
\end{figure}

In the \cirrs{} we use a rectangular contour $\mathcal {C}$ whose lower-left and upper-right vertices are $140$ and $335.4+50\im$, respectively.
The sampling points $z_i$ are set to be the nodes of the Gauss-Legendre quadrature rules on the four sides of $\mathcal {C}$. The \cirrs{} is ran for two cases with the same $L=4$ but different $N$, i.e., $N=30$ and $N=36$, respectively. In the case of $N=30$, the sampling points are distributed on the four sides in the ``10-5'' manner: 10 points Gauss-Legendre rule is used on each of the two longer sides, and 5 points rule is used on each of the two shorter sides. In the case of $N=36$, the sampling points are distributed in the ``12-6'' manner.
It turns out that there are 22 eigenvalues within this contour, including all the 21 eigenvalues of the NLEIGS. Figure \ref{fig-gun-ecs} shows the contour, the distributions of the 22 eigenvalues and the sampling points of the case $N=30$. Note that the eigenvalue marked by the circle is not included in the upper half-disk of the NLEIGS method.

\begin{figure}
\centering
\begin{minipage}[t]{.54\textwidth}
\centering
\vspace{0pt}
\epsfig{figure=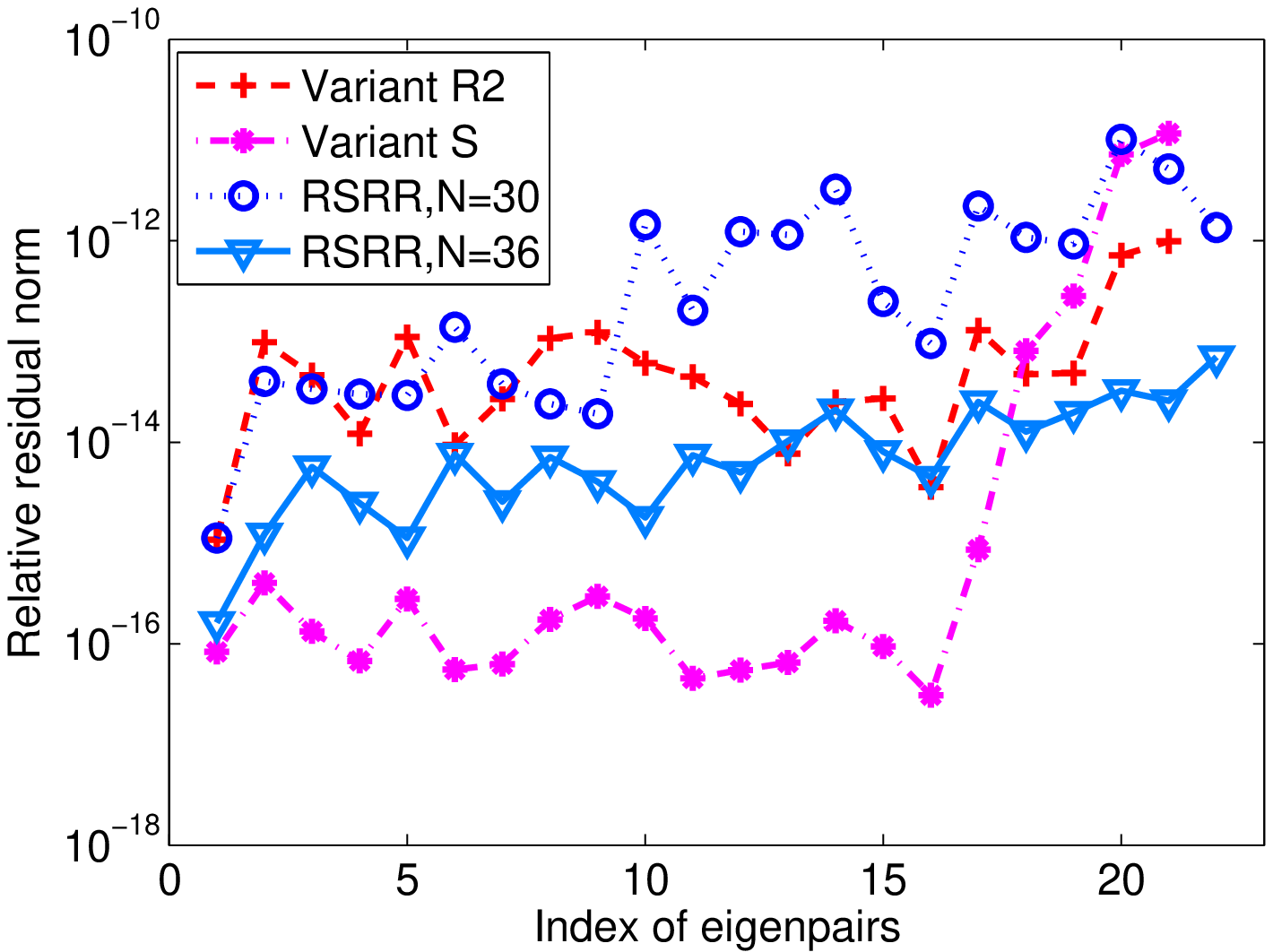,width=0.99\textwidth}
\end{minipage}
\hfill
\begin{minipage}[t]{.44\textwidth}
\centering
\vspace{20pt}
   \begin{tabular}[bht]{ccc}\hline
      Algorithms  & CPU time & Memory usage \\ \hline
      Variant R2  & 35.6 s     & $\thicksim$450 MB\\
      Variant S   & 13.1 s    & $\thicksim$435 MB\\
      \cirrs{}, $N=30$    & 29.1 s    & $\thicksim$140 MB\\
      \cirrs{}, $N=36$    & 34.7 s    & $\thicksim$140 MB\\ \hline
    \end{tabular}
\end{minipage}
\caption{The relative residual norms of the computed eigenpairs, total CPU times and memory usage of the NLEIGS and \cirrs{} in the \textit{`gun'} example.}
\label{fig-gun-errs}
\end{figure}

The relative residual norm (see \cite{guttel2014nleigs} for the definition) of the computed eigenpairs, total CPU times and memory usage of the NLEIGS and \cirrs{} are shown in Figure \ref{fig-gun-errs}.
First, the residuals of the NLEIGS and \cirrs{} are all below $10^{-10}$, and the accuracy of the \cirrs{} can be considerably improved by increasing $N$ from $30$ to $36$.
Second, the total CPU time of the \cirrs{} with $N=36$ is close to that of the \emph{Variant R2}, because the similar numbers of matrix inversions are required in the two methods. Third, the memory usage of the \cirrs{} is much lower than the NLEIGS, since in the \cirrs{} the linear equations are solved independently and no LU factorizations need to be saved for reuse, and moreover, only $N\cdot L$ solution vectors of dimension $n$ have to be saved, while in the NLEIGS the solution vectors would be of a much larger dimension.

From this example, we see that the \cirrs{} is comparable to the \emph{Variant R2} of the NLEIGS in terms of the accuracy and CPU time. In fact, these two methods have a close connection: both of them constructing the eigenspaces using the probed resolvent $T(z_i)^{-1}U$ on a series of sampling points $z_i$; see the discussion at the end of Section \ref{S-eigenspace-samp}.
The \emph{Variant S} of the NLEIGS is the fastest approach for this specific example, but we would like to remark that the performance of the \emph{Variant S} depends heavily on the number and distribution of the shifts in the rational Krylov method, and there is an example in \cite{guttel2014nleigs} showing that the \emph{Variant S} can be inferior to the \emph{Variant R2}. Finally, we stress again that the \cirrs{} is more suitable for parallelizations.

\subsection{BEM example}\label{S-ne-bem}

\begin{figure}[htb]
\centering
  \epsfig{figure=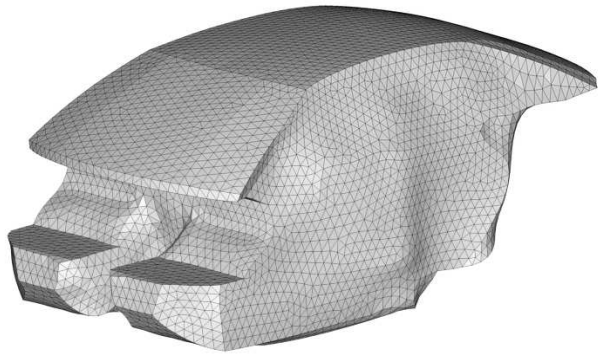,width=0.5\textwidth}
\caption{BEM mesh of the car cabin model}
\label{fig-car-mesh}
\end{figure}

To demonstrate the potential of the \cirrs{} method in BEM applications, a more engineering oriented example, the acoustic modal analysis of a car cabin cavity, is performed.
The boundary of the car cabin cavity is partitioned into 9850 triangular quadratic elements, leading to a BEM model with 59100 DOFS; see Figure \ref{fig-car-mesh}. The entire boundary is assumed to be rigid, i.e., $q = 0$. We solve for the eigenvalues in the interval $[40,500]$Hz. The contour is chosen to be the ellipse defined by $\gamma = 270, \, a = 230$ and $b = 0.1a$. Other controlling parameters of the \cirrs{} are set as: $L = 2$, $N = 200$. It turns out that there are 64 eigenvalues in this interval.

The BEM linear systems $T(z_i)^{-1} U$ are solved by using the GMRES solver with the ILU preconditioner. Both the accuracy of the fast BEM and the convergence tolerance of the GMRES solver are set to be $10^{-6}$. The matrix $T(\lambda)$ is interpolated by using Chebyshev polynomials of degree $d = 40$. In solving the reduced NEPs, parameters $N_S = 500$ and $K = 2$ are used.

%\subsubsection{Car cabin}

\begin{figure}[htb]
\centering
\begin{subfigure}[t]{.5\textwidth}
  \centering
  \epsfig{figure=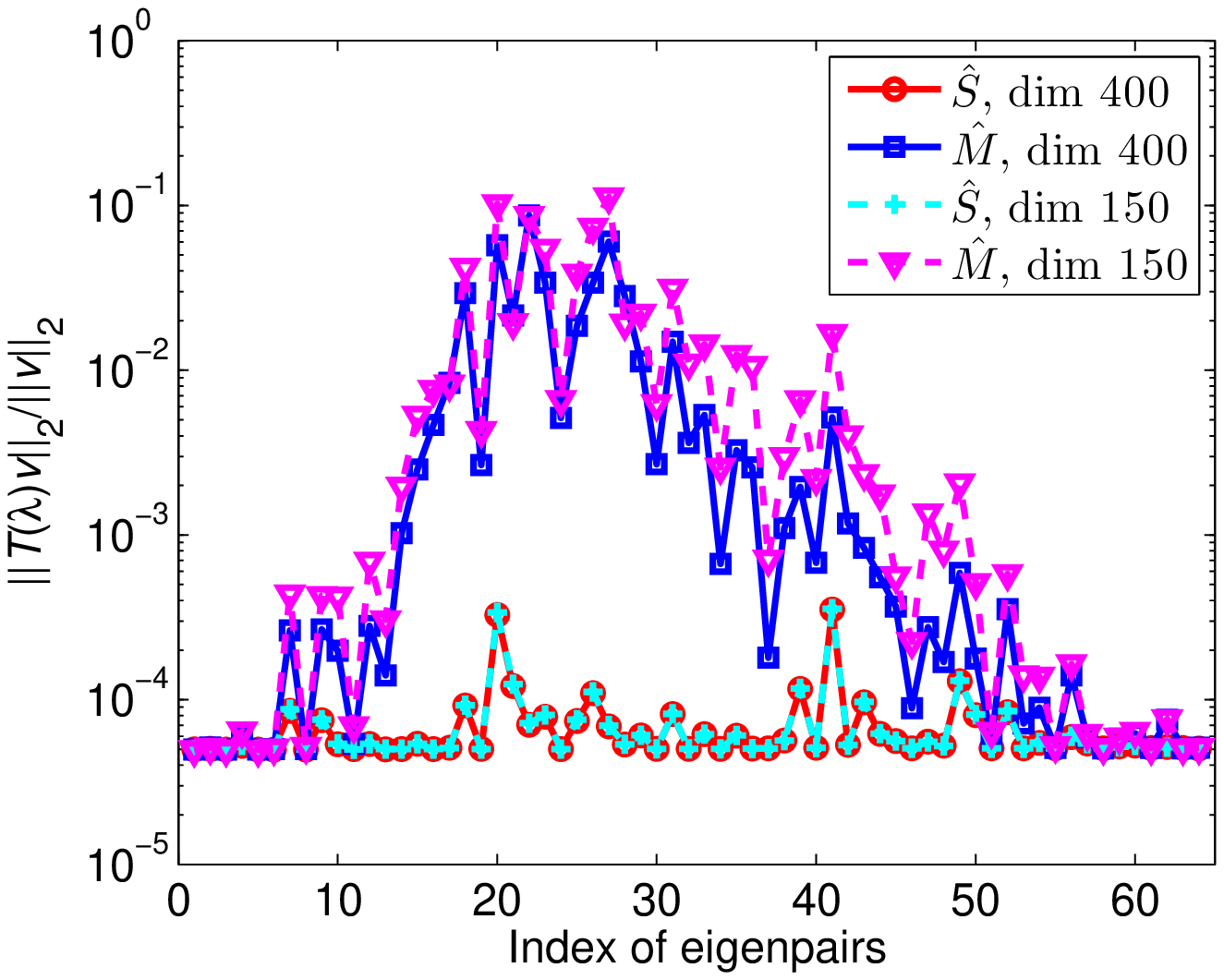,width=0.975\textwidth}
  \caption{Residuals of eigenpairs}
  \label{fig-car-err}
\end{subfigure}%
\begin{subfigure}[t]{.5\textwidth}
  \centering
  \epsfig{figure=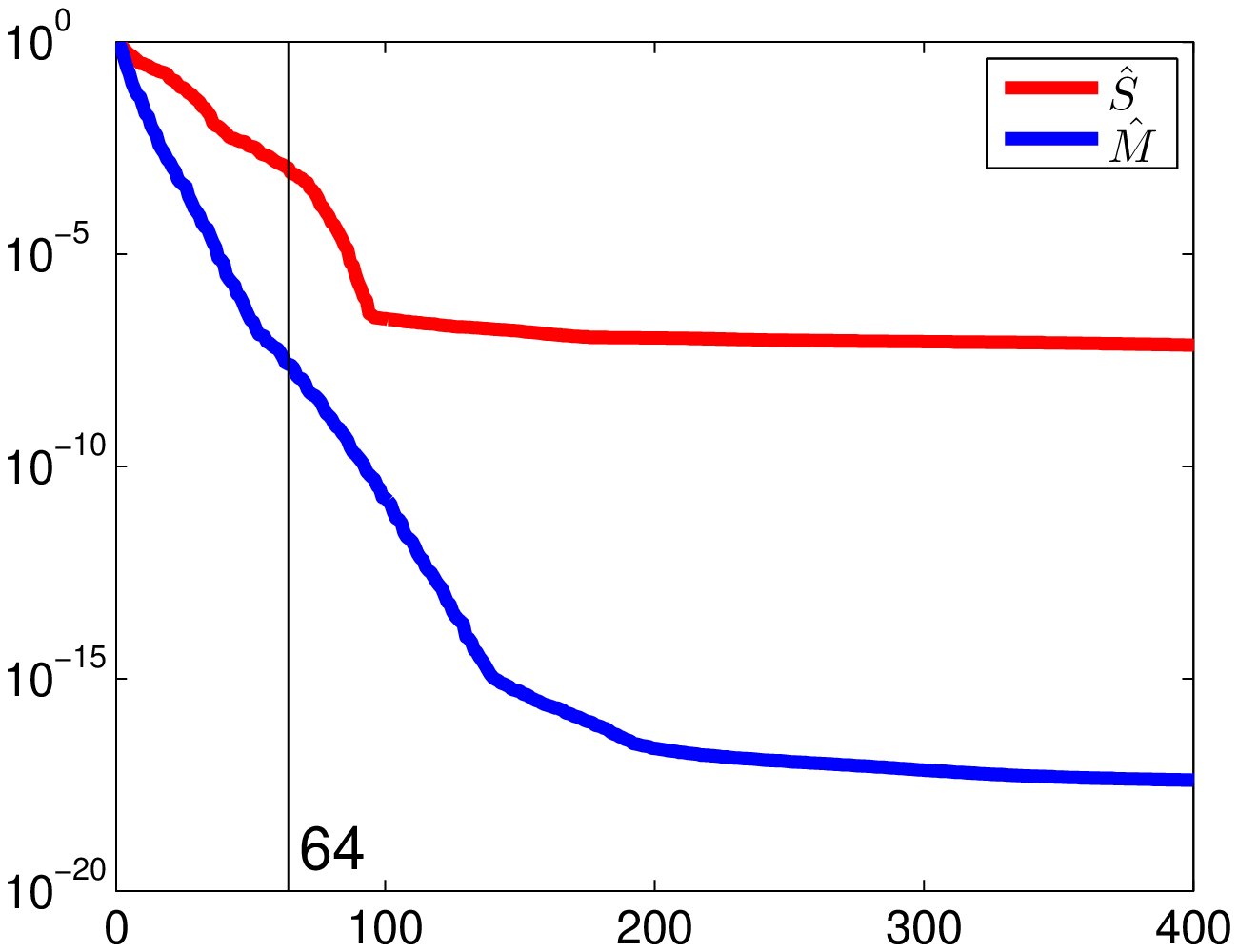,width=0.975\textwidth}
  \caption{Singular values of $\hat S$ and $\hat M$}
  \label{fig-car-svs}
\end{subfigure}
\caption{Comparisons of the sampling scheme and moment scheme for constructing eigenspaces for the car cabin model.}
\label{fig-car-errsvs}
\end{figure}

 To compare the performance of the \cirrs{} and \cirrm{} in the BEM. We consider two cases. In the first case the SVDs of $\hat S$ and $\hat M$ are both not truncated, so the dimensions of the eigenspaces are all 400, while in the second case the SVDs of  $\hat S$ and $\hat M$ are truncated so that the dimensions of the eigenspaces are all 150. The relative residuals of the computed 64 eigenpairs are shown in Figure \ref{fig-car-err}. Similar to the previous results, the sampling scheme performs much better than the moment scheme. When decreasing the dimension of the eigenspace from 400 to 150, no obvious increase of error is observed for the sampling scheme, indicating that the eigenspace information is effectively concentrated to the singular vectors of $\hat S$ associated with the largest singular values. But for the moment scheme, increasing the dimension of the eigenspace slightly improved the results; however, even when using the largest eigenspace (dimension 400), most eigenpairs at the middle part of the spectral are still of very low accuracy.
Figure \ref{fig-car-svs} shows the behaviors of the singular values of $\hat S$ and $\hat M$.

\subsection{FEM examples} \label{S-ne-fem}

The \cirrs{} algorithm is applied to solve a typical NEP in the engineering finite element analysis. We consider two models. The first one is a viscoelastic sandwich plate model, mainly used to confirm the superiority of the newly proposed sampling scheme in eigenspace generation. The second one is a large-scale model consisting of around one million DOFS.
The finite element analyses in these two examples are conducted using \textsc{Ansys}\textregistered{}. The linear systems in computing $T(z_i)^{-1} U$ are solved by using the sparse LU decomposition solver in \textsc{Ansys}\textregistered{}, and the accuracy of the solutions is always above 6 significant digits. In solving the reduced NEPs we use $N_S = 500$, $K = 2$.

\subsubsection{Viscoelastic sandwich plate}

\begin{figure}[htb]
\setlength{\unitlength}{1mm}
\centering
\begin{overpic}[width=0.5\textwidth]%
{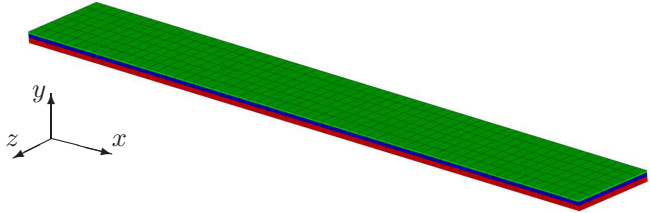}
\begin{picture}(1,1)
\put(1,10){\vector(4,-1){8}}
\put(1,10){\vector(0,1){6}}
\put(1,10){\vector(-2,-1){5}}
\put(9,9){$x$}
\put(-1.5,15.5){$y$}
\put(-5,9){$z$}
\end{picture}
\end{overpic}
\caption{The FEM model of the viscoelastic sandwich plate. The model consists of three layers from top to bottom: the constraining layer, the damping layer and the base layer.}
\label{fig-beam-mesh}
\end{figure}

The viscoelastic sandwich plate model shown in Figure \ref{fig-beam-mesh} is selected from \cite{LHW14} (Section 4.1). It consists of three layers. The base layer and constraining layer are Steel and Aluminum, respectively. In between these two layers is a layer of viscoelastic damping material which is modeled by using the Biot model. The plate is fixed at one side in the longest direction and free at the other sides. The sizes and material parameters are given in Table \ref{tab-beam-paras}.
The structure is discretized by using SOLID186 element (20-node hexahedral element exhibiting quadratic displacement behavior), leading to an analysis model with 18180 DOFS. The system matrix of the NEP is given by
\begin{equation}\label{eq-ne-dampingmodel}
T(z) = z^2 M + G(z)K_{\rm{v}} + K_{\rm{s}}, \quad \mathrm{and} \quad G(z) = G_\infty  \sum_{k=1}^3 {a_k z \over z + b_k} ,
\end{equation}
where, $M$ and $K_{\rm{s}}$ are the conventional mass matrix and stiffness matrix, $K_{\rm{v}}$ is the unit viscoelastic stiffness matrix, the three matrices are all of the dimension 18180, $G_\infty = 3.441 \times 10^{5}\,\rm{Pa}$, $a_1 = 2.06$, $a_2 = 67.1985$, $a_3 = 506.9457$, $b_1 = 193.39$ rad$/$s, $b_2 = 16345$ rad$/$s and $b_3 = 485918.4$ rad$/$s.

\begin{table}[!h]
\begin{center}
\caption{The geometry and material parameters of the three-layer sandwich plate \cite{LHW14}.}\label{tab-beam-paras}
\vspace{-.8\baselineskip}
\begin{tabular*}{0.8\textwidth}{@{\extracolsep{\fill}}llll@{}}\toprule
Parameters            &  Base layer          & Damping layer          & Constraining layer    \\
\hline
$x$-length (mm)       & 400                   & 400                   &  400  \\
$y$-length (mm)       & 3                     & 2                     &  1.2 \\
$z$-length (mm)       & 50                    & 50                    &  50  \\
Elastic modulus (Pa)  & $2.085\times 10^{11}$ & $1.0254\times 10^{6}$ &  $6.76\times 10^{10}$ \\
Density (kg$/$m$^3$)  & 7840                  & 789.5                 &  2803 \\
Poisson's ratio         & 0.3                   & 0.49                  &  0.3 \\
\bottomrule
\end{tabular*}
\end{center}
\end{table}

\begin{figure}[htb]
\centering
\epsfig{figure=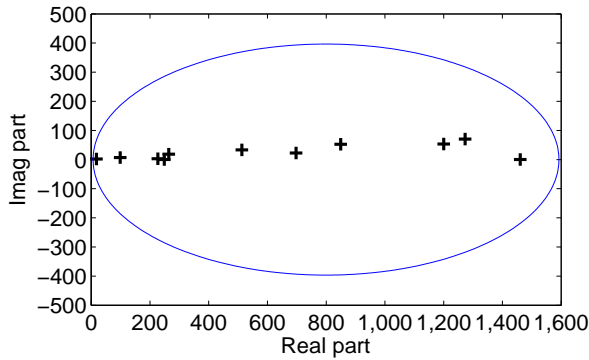,width=0.5\textwidth}
\caption{Computed eigenvalues ($\scriptstyle\pmb{+}$) and the contour ($-$) of the viscoelastic sandwich plate model}
\label{fig-beam-evs}
\end{figure}

\begin{figure}[htb]
\centering
\begin{subfigure}[t]{.5\textwidth}
  \centering
  \epsfig{figure=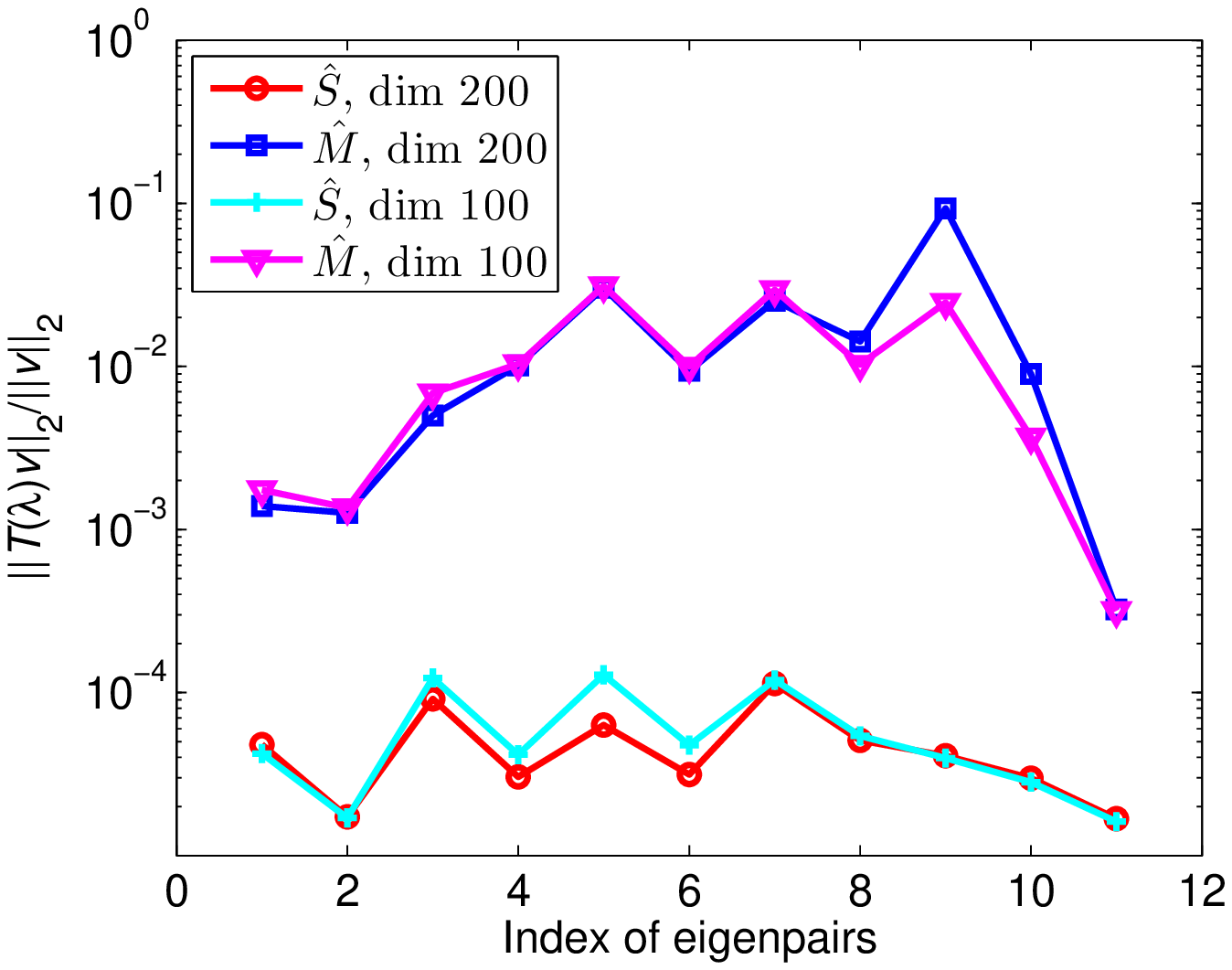,width=0.975\textwidth}
  \caption{Residuals of eigenpairs }
  \label{fig-beam-err}
\end{subfigure}%
\begin{subfigure}[t]{.5\textwidth}
  \centering
  \epsfig{figure=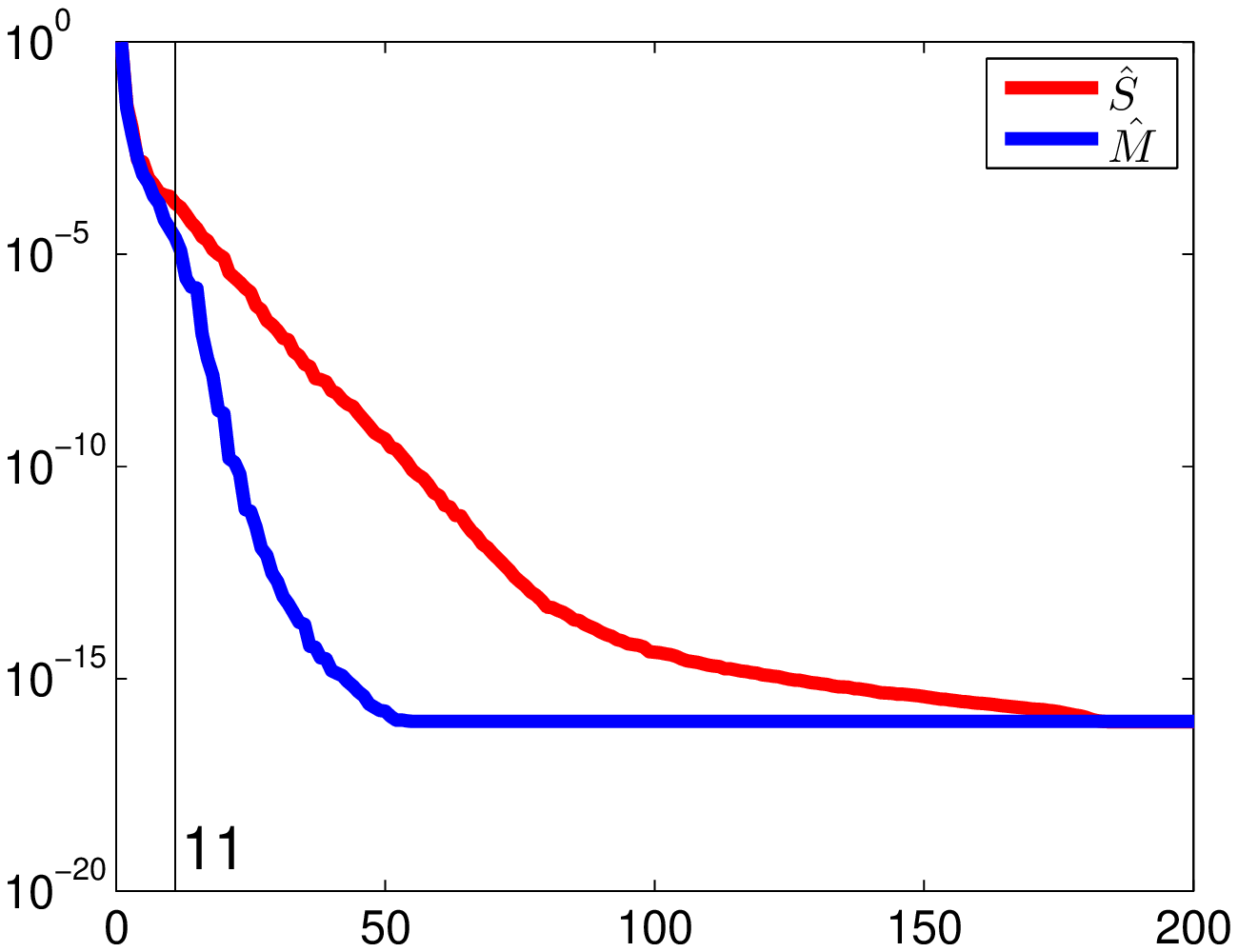,width=0.975\textwidth}
  \caption{Singular values of $\hat S$ and $\hat M$}
  \label{fig-beam-svs}
\end{subfigure}
\caption{Comparison of the sampling scheme and moment scheme for the viscoelastic sandwich plate model. In the left plot, ``dim'' indicates the dimension of the eigenspace. }
\label{fig-beam-errsvs}
\end{figure}

We solve for the natural frequencies in the interval $[6.4, 1592.4]$Hz. The contour is chosen to be the ellipse defined by $\gamma = 799.4, \, a = 793,\, b = 0.5a$. Other controlling parameters of \cirrs{} and \cirrm{} are set as: $L = 2$, $N = 200$, $K'=N$. We obtain 11 eigenpairs by using the \cirrs{}; see Figure \ref{fig-beam-evs} for the distribution of the computed eigenvalues. The relative residual $||T(\lambda)v||_2/||v||_2$ of the computed eigenpairs is illustrated in Figure \ref{fig-beam-err}, and the behavior of the singular values of the matrices $\hat S$ and $\hat M$ is illustrated in Figure \ref{fig-beam-svs}. Note that both $\hat S$ and $\hat M$ should have 400 singular values, but only the first 200 singular values are shown in Figure \ref{fig-beam-svs} because the remainders are almost constant around $10^{-16}$.

To test the performance of $\hat S$ and $\hat M$ in generating eigenspace, the SVDs of them are truncated so that the dimensions of the resultant eigenspaces are 200 and 100, respectively.
From Figure \ref{fig-beam-err}, one can see that the residuals of the sampling scheme are typically below $10^{-4}$, one order lower than those of the moment scheme. In addition, for both the two schemes, the change in the dimension of the eigenspaces does not lead to obvious changes in the residuals of the eigenpairs. In view of the singular value behavior in Figure \ref{fig-beam-svs}, this indicates that it is sufficient to truncate the SVD by $\delta = 10^{-14}$. In a nutshell, the sampling scheme can generate much better eigenspaces than the moment scheme, provided that the SVD of $\hat S$ is truncated by $\delta < 10^{-14}$. Our further numerical tests also show that the residuals of the sampling scheme increases for a larger $\delta > 10^{-14}$. Similar behavior can also be observed in the moment scheme.

\subsubsection{Payload attach fitting} \label{S-S-paf}

\begin{figure}[htb]
\setlength{\unitlength}{1mm}
\centering
\begin{overpic}[width=0.45\textwidth]%
{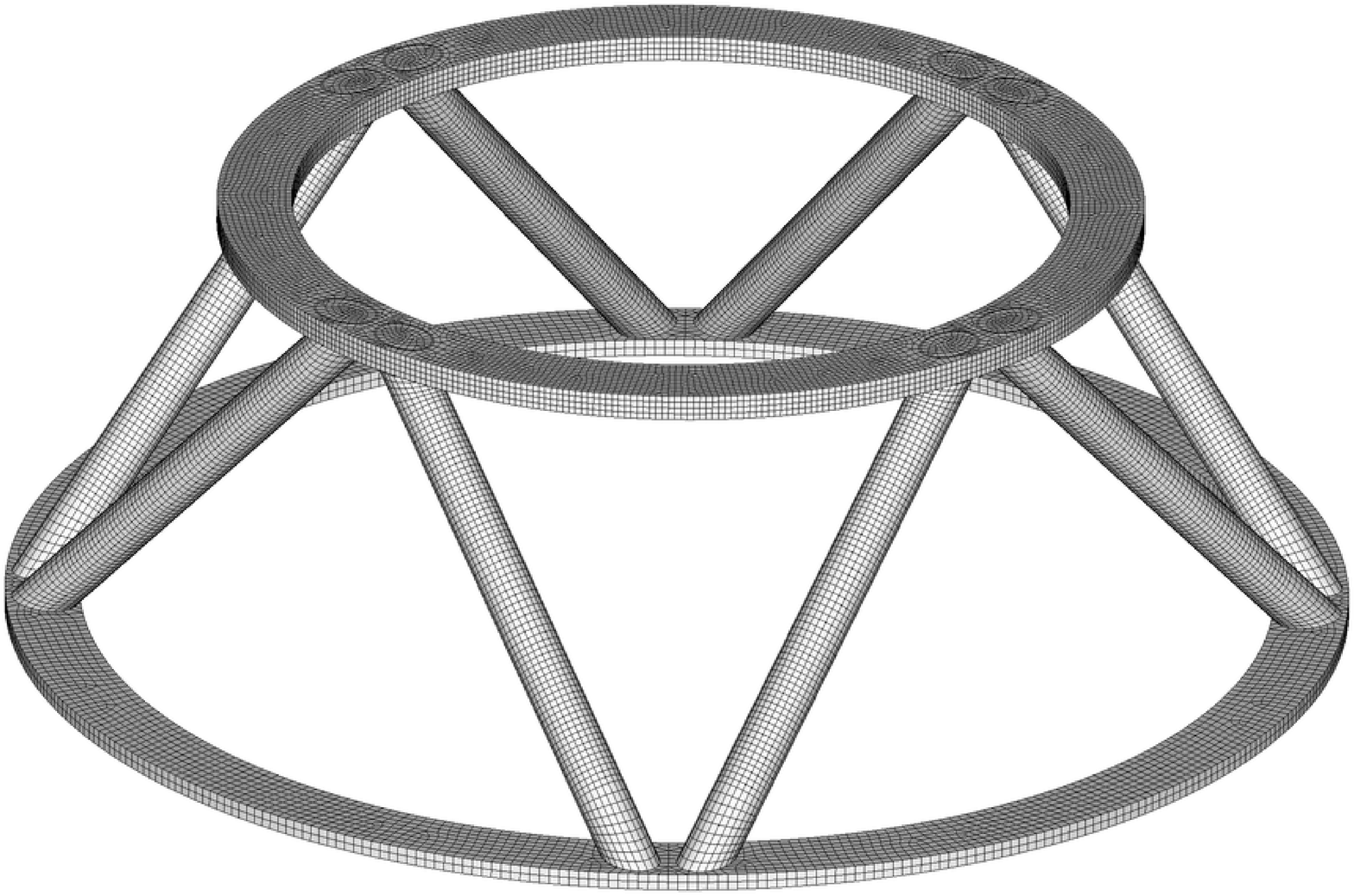}
\begin{picture}(1,1)
\put(34.5,17){\vector(1,0){6}}
\put(34.5,17){\vector(0,1){6}}
\put(34.5,17){\vector(-1,-1){3.6}}
\put(40,19){$y$}
\put(31.5,23){$z$}
\put(30,16){$x$}

\put(12,13){\vector(-2,-1){4}}
\put(6.5,14){{\footnotesize Lower Flange}}
\put(6,45){\vector(1,-1){4}}
\put(-4,47){{\footnotesize Upper Flange}}
\put(67,36){\vector(-1,-1){4}}
\put(62.8,40){{\footnotesize Damping}}
\put(63,37){{\footnotesize Cylinder}}
\end{picture}
\end{overpic}\hspace{0.5cm}
\epsfig{figure=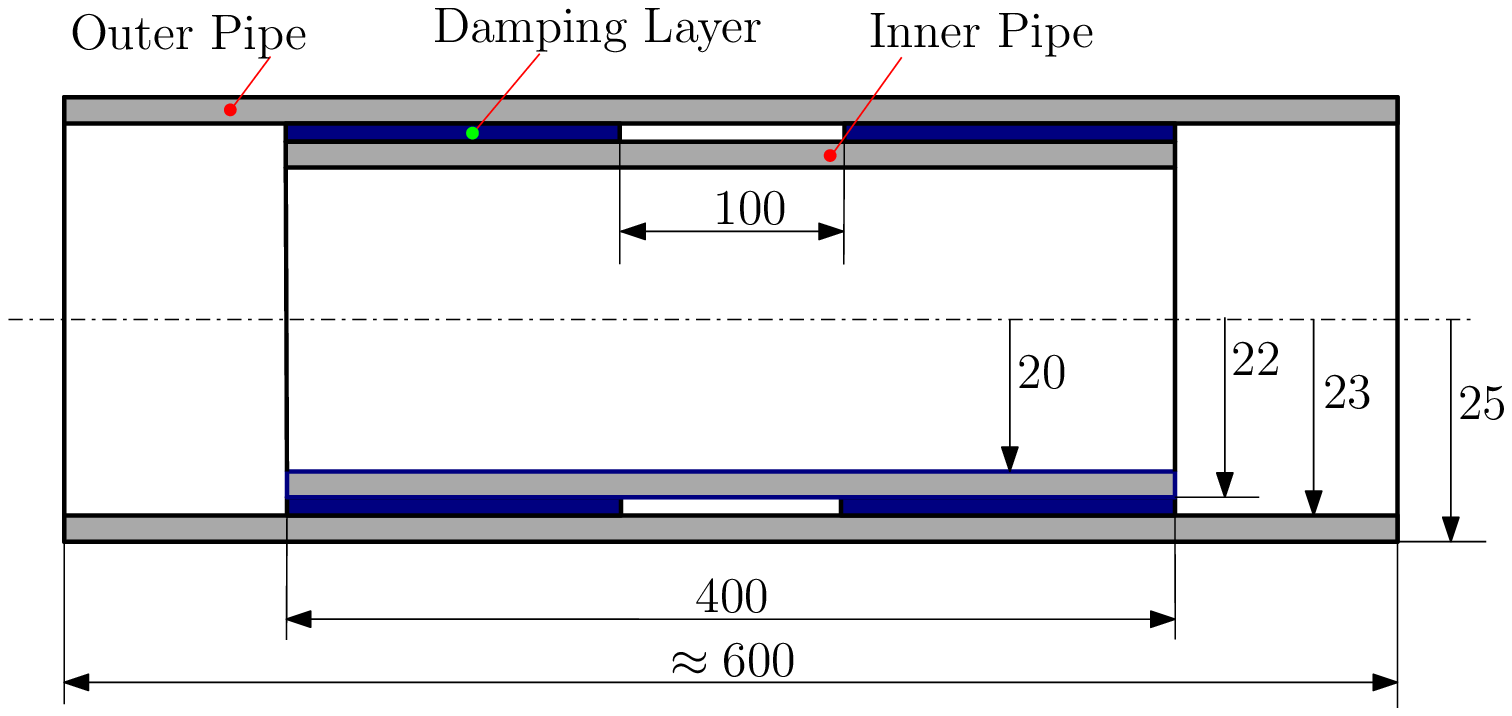,width=0.45\textwidth}
\caption{The payload attach fitting model (Unit: mm). The model consists of an upper flange, a lower flange and eight damping cylinders. The dimensions of the upper flange are: inner radius 365, outer radius 435, thickness 24. The dimensions of the lower flange are: inner radius 560, outer radius 630, thickness 16. The structure and dimensions of the damping cylinders are shown in the left diagram.}
\label{fig-paf-model}
\end{figure}

As a large-scale example, we consider the payload attach fitting model illustrated in Figure \ref{fig-paf-model}. This kind of structures usually serve as isolators
between the satellite and the launch vehicle to reduce the satellite vibration caused by the launch-induced dynamic loads. The structural material is Aluminum with elastic modulus $70\,\rm{GPa}$, density $2770\, \rm{kg}/\rm{m}^3$ and Poisson's ratio $0.3$. The viscoelastic damping material in the damping cylinders is ZN1 polymer which is modeled by the Biot model. The density and Poisson's ratio of the damping material are $970\, \rm{kg}/\rm{m}^3$ and $0.49$, respectively. The structure is fixed at the lower surface of the lower flange and all the other surfaces are free of traction. The model is discretized by using SOLID186 element, and the total number of DOFS is 1005648. The expression of this NEP is the same as \eqref{eq-ne-dampingmodel}, but with
$G(z) = G_{\infty} \left(1 + \sum_{k=1}^4 {a_k z \over z + b_k}\right)$, where $G_{\infty} = 362750\,\rm{Pa}$, $a_1 = 0.762063$, $a_2 = 1.814626$, $a_3 = 84.93828$, $a_4 = 4.869723$, $b_1 = 53.72964$, $b_2 = 504.5871$, $b_3 = 29695.64$ and $b_4 = 2478.43$.

We compute the natural frequencies below 640Hz. A rough modal analysis of the undamped structure shows that the lowest natural frequency is around 200 Hz. So we set the frequency interval of the \cirrs{} to be $[160, 640]$Hz. The ellipse contour with $\gamma = 400, \, a = 240$ and $b = 0.5a$ is used accordingly. Besides, we let $L = 2$ and $N = 100$. Thus there are totally 200 FEM linear systems to be solved.

\begin{figure}[htb]
\centering
\begin{subfigure}[t]{.32\textwidth}
  \centering
  \epsfig{figure=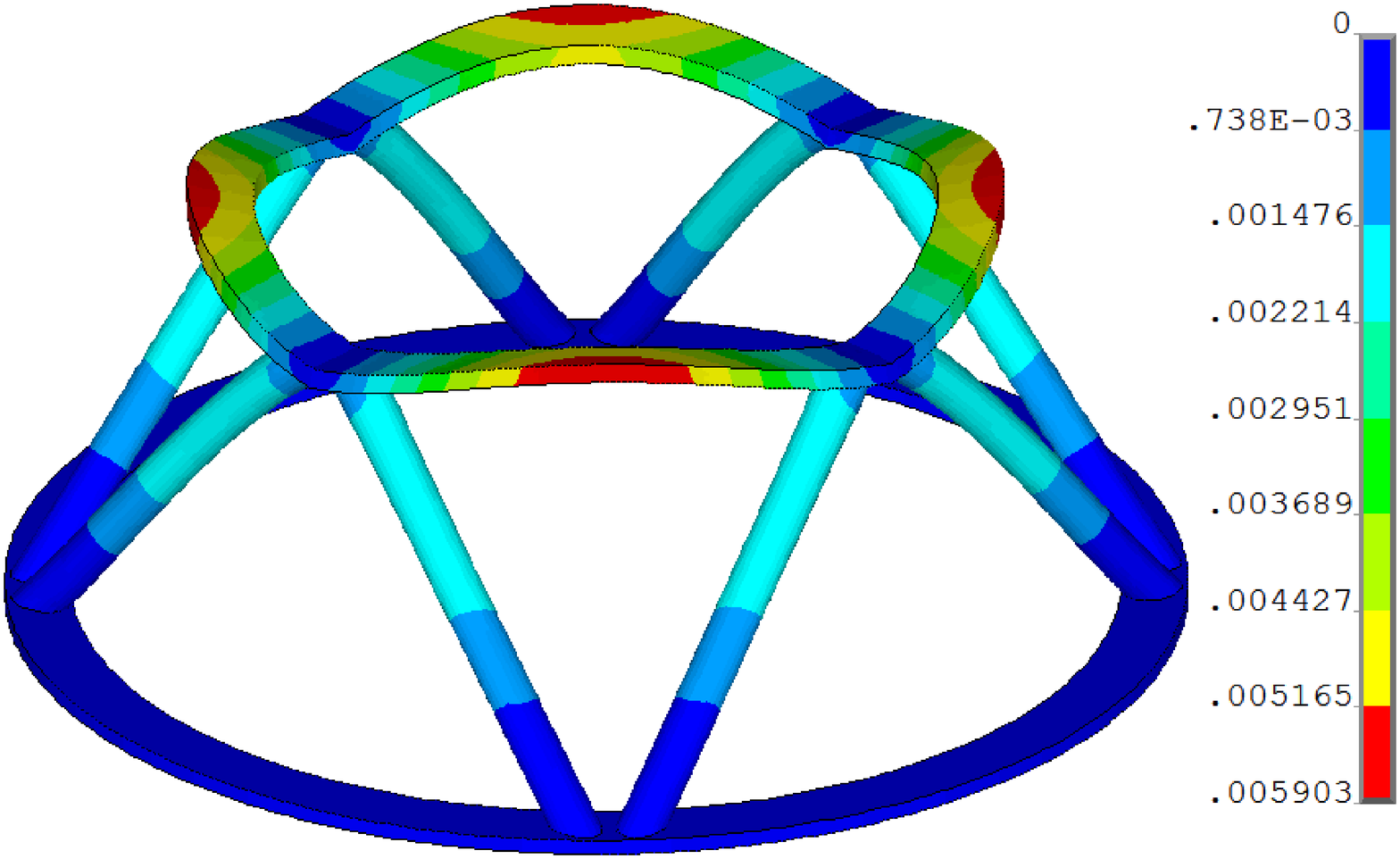,width=0.99\textwidth}
  \caption{Mode 1}
  \label{fig-paf-m1}
\end{subfigure}%
\begin{subfigure}[t]{.32\textwidth}
  \centering
  \epsfig{figure=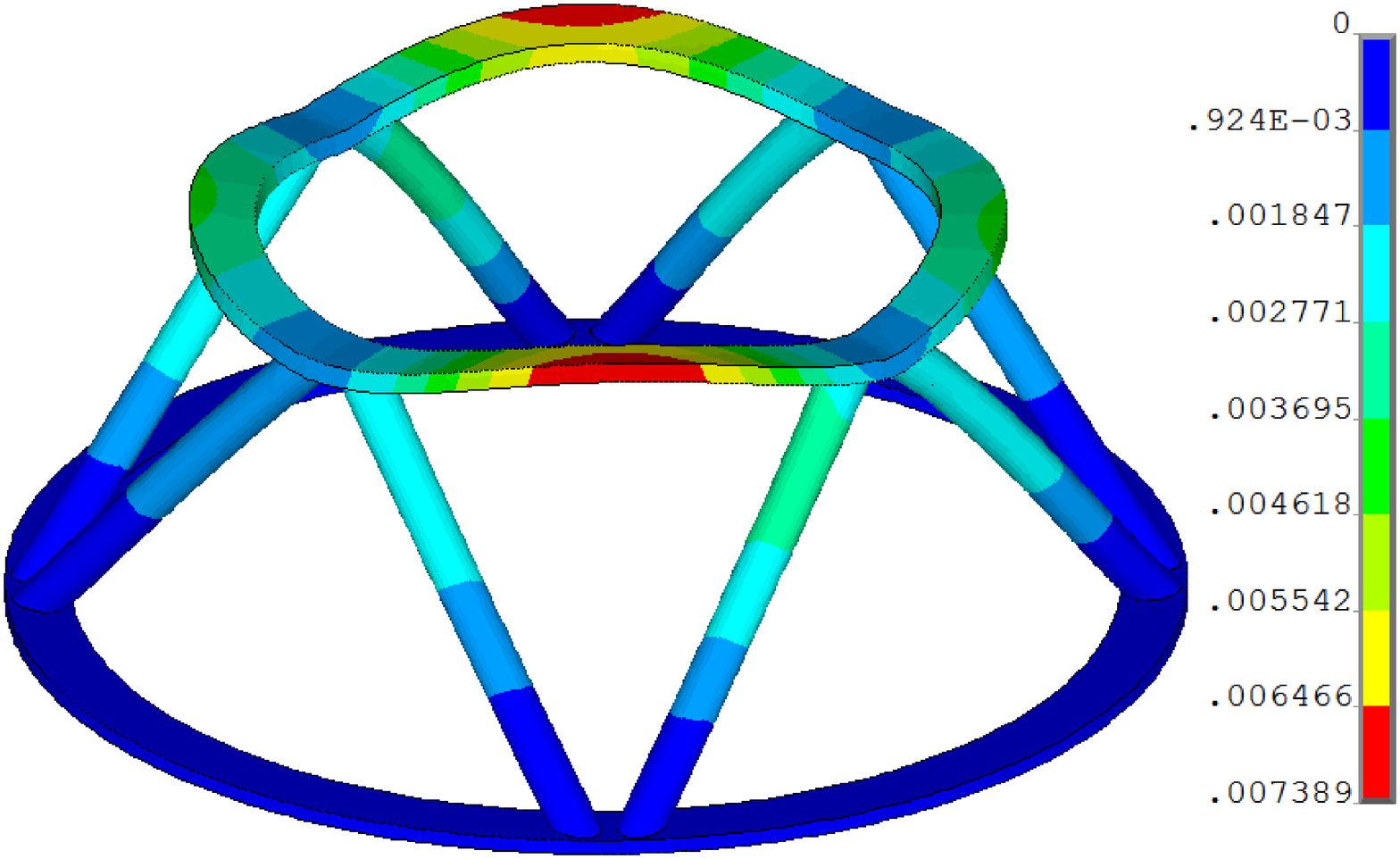,width=0.99\textwidth}
  \caption{Mode 2}
  \label{fig-paf-m2}
\end{subfigure}
\begin{subfigure}[t]{.32\textwidth}
  \centering
  \epsfig{figure=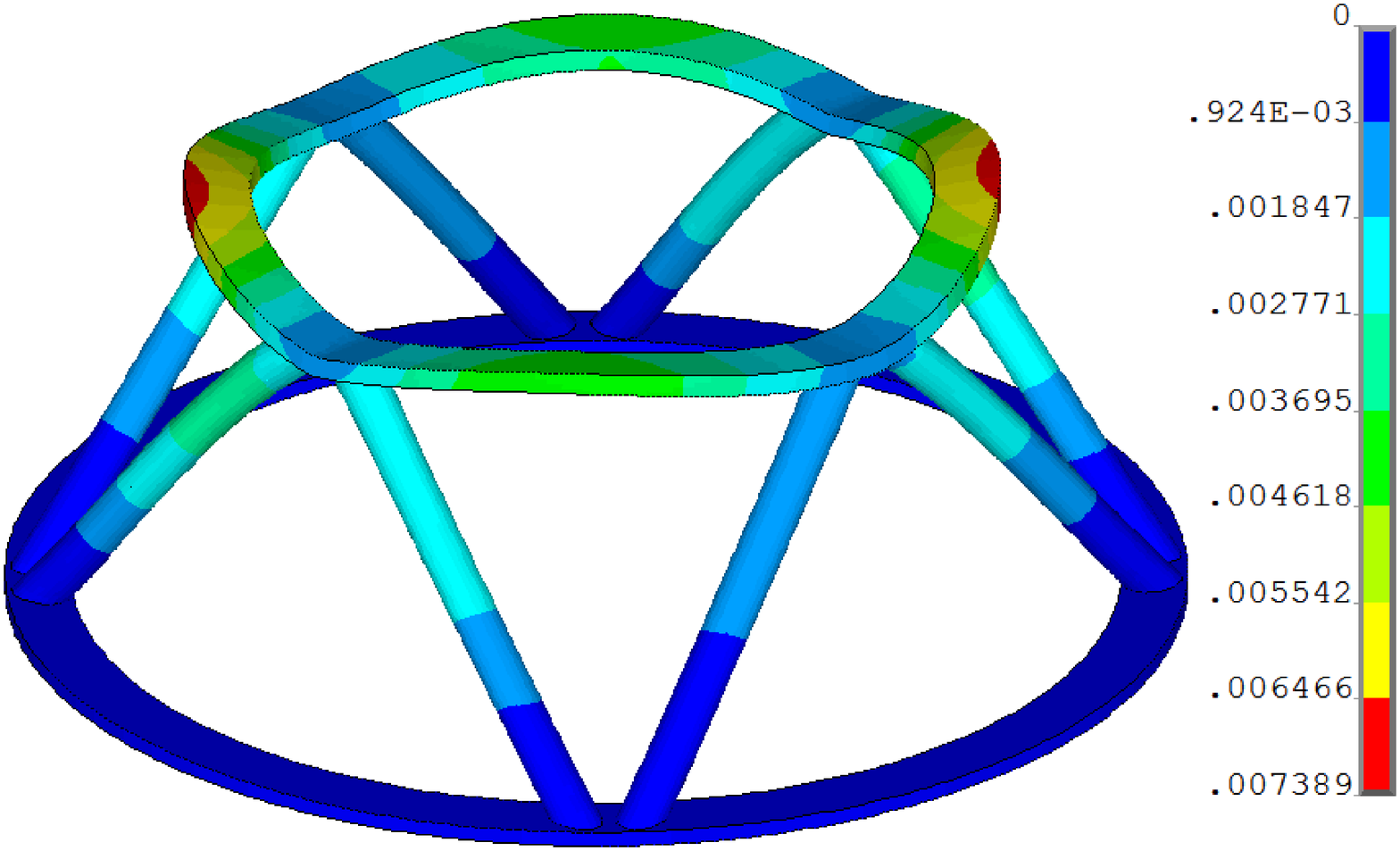,width=0.99\textwidth}
  \caption{Mode 3}
  \label{fig-paf-m3}
\end{subfigure}%

\begin{subfigure}[t]{.32\textwidth}
  \centering
  \epsfig{figure=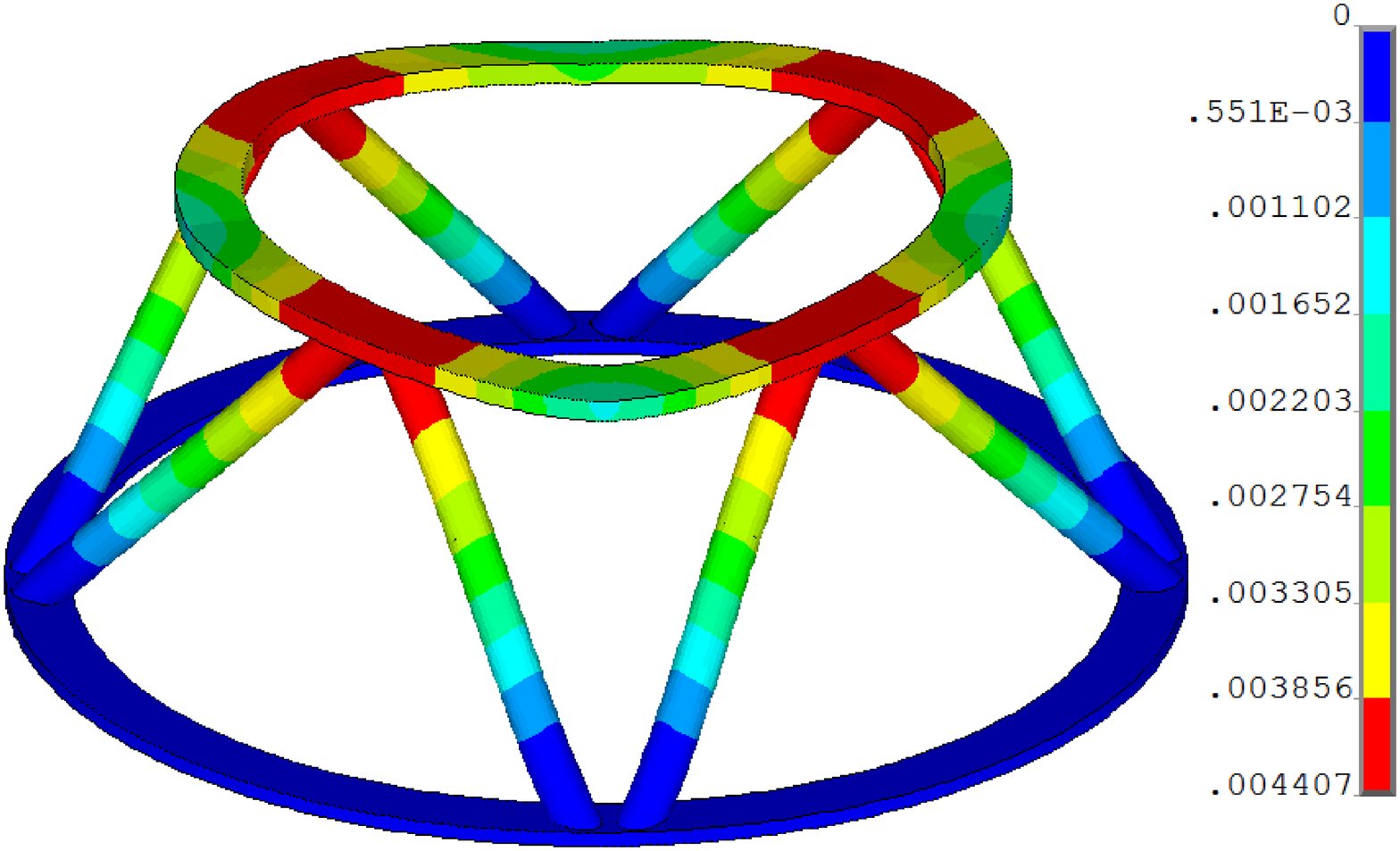,width=0.99\textwidth}
  \caption{Mode 4}
  \label{fig-paf-m4}
\end{subfigure}
\begin{subfigure}[t]{.32\textwidth}
  \centering
  \epsfig{figure=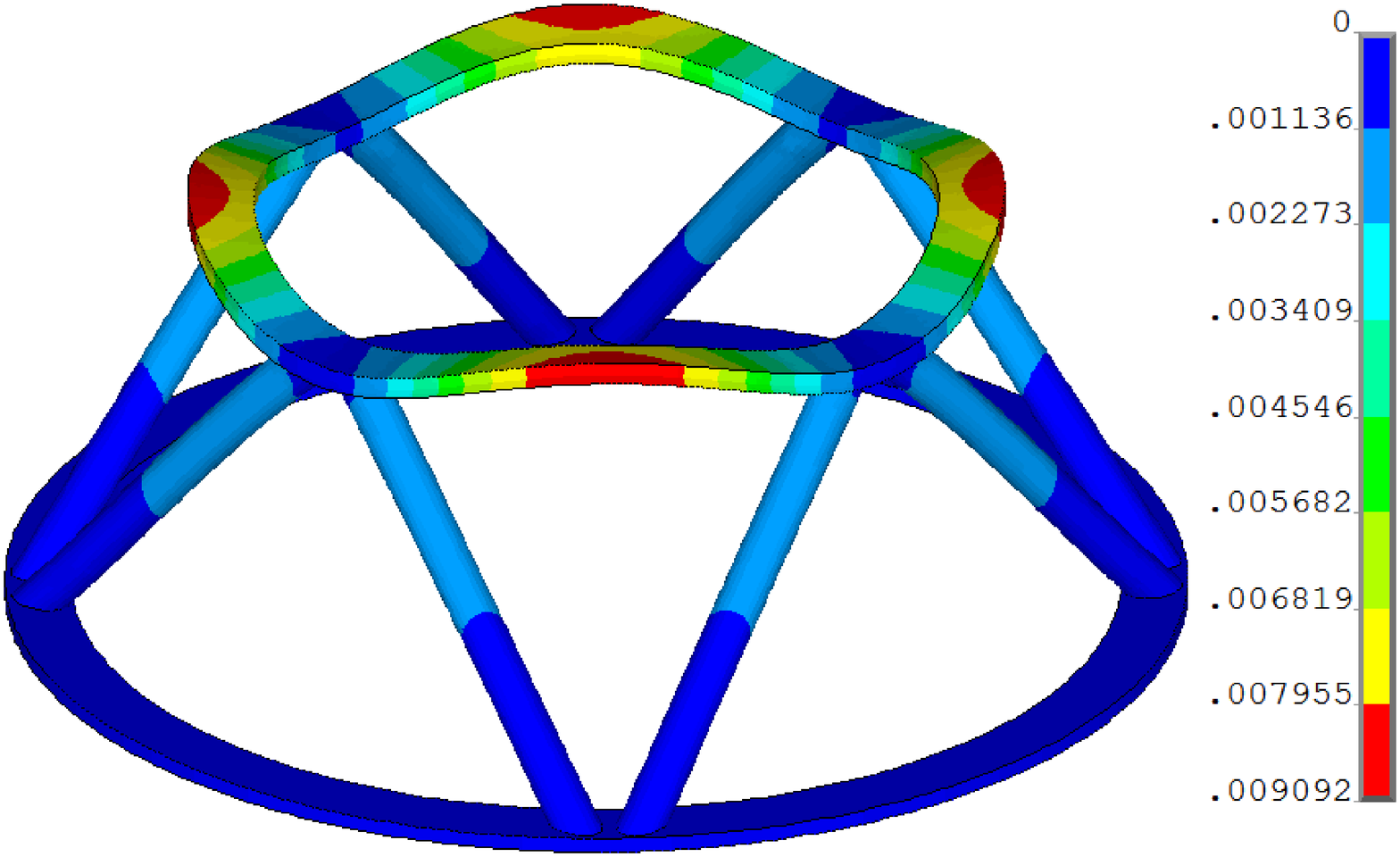,width=0.99\textwidth}
  \caption{Mode 5}
  \label{fig-paf-m5}
\end{subfigure}
\begin{subfigure}[t]{.32\textwidth}
  \centering
  \epsfig{figure=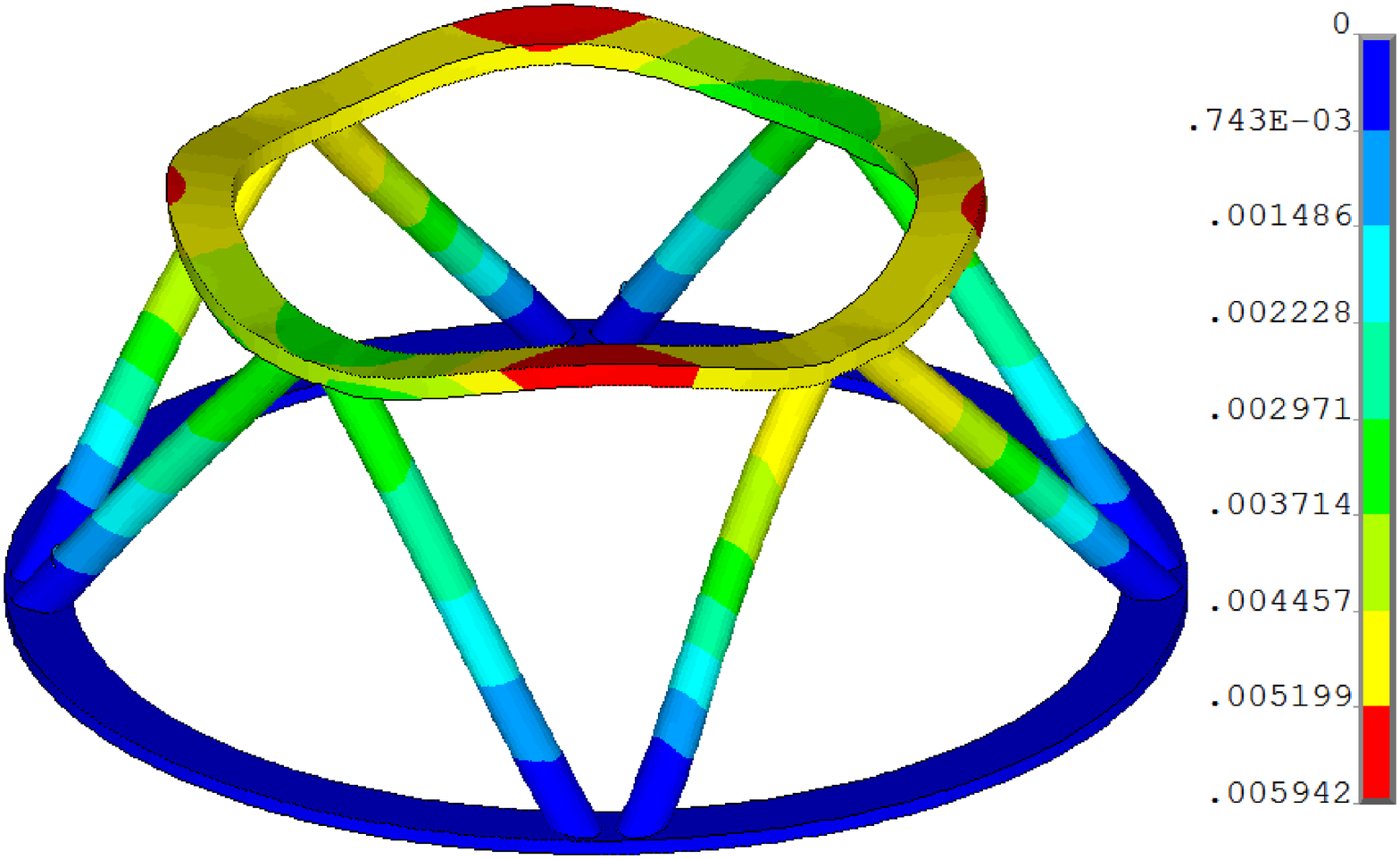,width=0.99\textwidth}
  \caption{Mode 6}
  \label{fig-paf-m6}
\end{subfigure}
\caption{The first 6 modes (amplitude) of the payload attach fitting model. }
\label{fig-paf-6modes}
\end{figure}

The computation was carried out on a Server with eight 8-core Intel Xeon (2.67GHz) processors. The whole solution process takes about 7 hours. There are 14 eigenpairs; see Table \ref{tab-paf-evs} for the eigenfrequencies and the corresponding residuals. Due to the symmetry of the structure, there are 3 couples of 2-fold eigenfrequencies. The largest residual is around $10^{-5}$, indicating a good accuracy of the \cirrs{} for solving real-world problems. The accuracy of the computed eigenpairs is mainly limited by that of the linear system solver (sparse LU decomposition) in \textsc{Ansys}\textregistered{}, because a largest relative error of around $10^{-6}$ was observed in solving the 200 linear systems. Nevertheless, an accuracy of 5 significant digits should be sufficient for most engineering applications.
Figure \ref{fig-paf-6modes} exhibits the amplitude of the first 6 modes of the model. Clearly, modes 2 and 3 correspond to a 2-fold frequency and so they are symmetric with each other. Modes 1, 4 and 5 correspond to simple eigenfrequencies, each of them is symmetric with respect to $z$-axis that passes through the centers of the lower and upper flanges. Analogously, mode 6 should be symmetric with mode 7 which is not illustrated.

\begin{table}[!h]
\begin{center}
\caption{The 14 computed eigenvalues and the residuals of the payload attach fitting model}\label{tab-paf-evs}
\vspace{-.8\baselineskip}
\begin{tabular*}{0.8\textwidth}{@{\extracolsep{\fill}}ccc||ccc@{}}\toprule
%&\multicolumn{3}{c}{rre}    &\multicolumn{2}{c}{Present}\\
$k$               &  $\lambda_k$    & $||T(\lambda)v||_2/||v||_2$ & $k$               &  $\lambda_k$    & $||T(\lambda)v||_2/||v||_2$    \\
\hline
1  & $206.88 + 0.99\im$   & $3.80\times 10^{-6}$  &  8  & $405.56 + 8.45\im$   & $1.93\times 10^{-6}$\\
2  & $228.06 + 1.79\im$   & $5.99\times 10^{-6}$  &  9  & $431.31 + 12.26\im$   & $3.60\times 10^{-6}$\\
3  & $228.06 + 1.79\im$   & $1.14\times 10^{-5}$ &  10 & $551.31 + 6.26\im$   & $4.80\times 10^{-7}$\\
4  & $251.95 + 0.16\im$   & $1.54\times 10^{-5}$ &  11 & $551.31 + 6.26 \im$   & $3.21\times 10^{-6}$\\
5  & $289.73 + 2.91\im$   & $1.70\times 10^{-7}$   &  12 & $573.46 + 2.98\im$   & $3.50\times 10^{-7}$\\
6  & $345.70 + 12.25\im$   & $2.93\times 10^{-6}$ &  13 & $607.30 + 15.48\im$   & $2.75\times 10^{-6}$\\
7  & $345.70 + 12.25\im$   & $3.46\times 10^{-6}$ &  14 & $630.17 + 17.92\im$   & $1.51\times 10^{-6}$\\
\bottomrule
\end{tabular*}
\end{center}
\end{table}

\section{Conclusions} \label{S-conclusion}

This paper is devoted to develop efficient numerical methods for solving large-scale NEPs in computational science and engineering. The outcome is a robust and reliable method, \cirrs{}, which is suitable for general NEPs of the form $T(\lambda)v = 0$, no matter the matrix $T(z)$ is sparse or dense, structured or unstructured. The method is based on the Rayleigh-Ritz procedure. We have proposed a new resolvent sampling scheme for generating reliable eigenspaces, and a modified block SS algorithm for the robust and accurate solution of the reduced NEPs. The main computational work of the \cirrs{} lies in the solution of independent linear systems with $T(z)$ on the sampling points, and in many areas fast linear solvers are available.

As applications, the \cirrs{} algorithm has been used to solve representative NEPs in engineering finite element and boundary element analyses. Note that the application to the latter is not straightforward. The Chebyshev interpolation technique has been employed to obtain approximate expressions for the matrix $T(z)$ and its reduced counterpart. The good performance of the newly proposed sampling scheme in generating eigenspaces and finally the robustness, reliability and accuracy of the \cirrs{} have been attested by several NEPs in a variety of applications. In addition, the potential of the \cirrs{} in large-scale simulation has been demonstrated by solving a FEM NEP with 1 million DOFS and a BEM NEP with 60 thousand DOFS. Besides the robustness and high accuracy, the \cirrs{} algorithm is very suitable for parallelization, and can be easily implemented into other programs and software.

\section*{Acknowledgements}

JX gratefully acknowledges the financial supports from the National Science
Foundations of China under Grants 11102154 and 11472217, Fundamental Research Funds for the Central Universities and the Alexander von Humboldt Foundation (AvH) to support his fellowship research at the Chair of Structural Mechanics, University of Siegen, Germany.
The authors wish to thank the valuable
comments of the anonymous reviewers which significantly improved the quality
of this manuscript.

\section*{References}

\bibliography{nep_rr_general_bib}

\end{document}